\documentclass[11pt]{amsart}
\usepackage{amsmath,amssymb,hyperref}
\usepackage{amsfonts}
\usepackage{enumerate}
\usepackage{amsmath, amsthm}
\usepackage{amscd}
\input xy
\xyoption{all}
\hypersetup{pdfpagemode=FullScreen, colorlinks=true}

\newtheorem{thm}{Theorem}[section]
\newtheorem{sthm}{Theorem}

\newtheorem{prop}[thm]{Proposition}
\newtheorem{lemma}[thm]{Lemma}

\newtheorem{cor}[thm]{Corollary}
\newtheorem{defi}[thm]{Definition}
\newtheorem{nota}[thm]{Notation}
\newtheorem{rem}[thm]{Remark}
\newtheorem{ex}[thm]{Example}
\newtheorem{qu}[thm]{Question}

\newcommand{\bc}{{\mathbb C}}
\newcommand{\bh}{{\mathbb H}}
\newcommand{\br}{{\mathbb R}}
\newcommand{\bq}{{\mathbb Q}}
\newcommand{\bz}{{\mathbb Z}}
\newcommand{\ra}{\rightarrow}
\newenvironment{pf}{\begin{trivlist}\item[]{\bf Proof:\ }}
{\mbox{}\hfill\rule{.08in}{.08in}\end{trivlist}}

\begin{document}
\title[Flexibility of surface groups in classical groups]
{Flexibility of surface groups in classical simple Lie groups}
\author{Inkang Kim and Pierre Pansu}
\date{}

\begin{abstract}
We show that a surface group of high genus contained in a classical
simple Lie group can be deformed  to become Zariski dense, unless
the Lie group is $SU(p,q)$ (resp. $SO^* (2n)$, $n$ odd) and the
surface group is maximal in some $S(U(p,p)\times U(q-p))\subset
SU(p,q)$ (resp. $SO^* (2n-2)\times SO(2)\subset SO^* (2n)$). This is
a converse, for classical groups, to a rigidity result of S.
Bradlow, O. Garc\'{\i}a-Prada and P. Gothen.
\end{abstract}

\footnotetext[1]{I. Kim gratefully acknowledges the partial support
of NRF grant  ((R01-2008-000-10052-0) and a warm support of IHES
during his stay.}

\maketitle

\hfill{\em Dedicated to Lionel B\'erard-Bergery, for his 65th birthday}

\section{Introduction}

Free groups are obviously flexible. In particular, a generic free subgroup in a real algebraic group is Zariski dense. This already fails for surface groups, although they are very flexible from other points of view. The first evidence came from the following result by D. Toledo.

\begin{thm}
\label{Toledo}
{\em (D. Toledo, 1979, 1989, \cite{To}).} Let $\Gamma$ be a discrete cocompact subgroup of $SU(1,1)$. Map $SU(1,1)$ as a $2\times 2$ block in $SU(1,n)$, $n\geq 2$. Then every neighboring homomorphism $\Gamma\to SU(1,n)$ is contained in a conjugate of $S(U(1,1)\times U(n-1))$.
\end{thm}

In fact, Toledo obtained a stronger, global result: a characterization of surface subgroups of $S(U(1,1)\times U(n-1))$ among surface subgroups of $SU(1,n)$ by the value of a characteristic class known as Toledo's invariant, which we now define.

Let $X$ be a Hermitian symmetric space, with K\"ah\-ler  form
$\Omega$ (the metric is normalized so that the minimal sectional
curvature equals $-1$). Let $\Sigma$ be a closed surface of negative
Euler characteristic, let $\Gamma=\pi_1 (\Sigma)$ act isometrically
on $X$. Pick a smooth equivariant map $\tilde{f}:\tilde{\Sigma}\to
X$.

\begin{defi}
Define the Toledo invariant of the action $\rho:\Gamma\to Isom(X)$ by
\begin{eqnarray*}
T_{\rho}=\frac{1}{2\pi}\int_{\Sigma}\tilde{f}^{*}\Omega.
\end{eqnarray*}
\end{defi}

Then
\begin{enumerate}
  \item $T_{\rho}$ depends continuously on $\rho$.
  \item There exists $\ell_X \in\bq$ such that $T_{\rho}\in\ell_X \bz$.
  \item $|T_{\rho}|\leq|\chi(\Sigma)|\mathrm{rank}(X)$.
\end{enumerate}
Inequality (3), known as the Milnor-Wood inequality for actions on Hermitian symmetric spaces, is due to J. Milnor, \cite{Milnor}, V. Turaev, \cite{Tu}, A. Domic and D. Toledo, \cite{DT}, J.-L. Clerc and B. \O rsted, \cite{cor}.

\begin{defi}
Actions $\rho$ such that $|T_{\rho}|=|\chi(\Sigma)|\mathrm{rank}(X)$ are called {\em maximal representations}.
\end{defi}

The following result generalizes Theorem \ref{Toledo}.

\begin{thm}
\label{BGPG}
{\em (L. Hern\`andez Lamoneda, \cite{HL}, S. Bradlow, O. Garc\'{\i}a-Prada, P. Gothen, \cite{BGG}, \cite{BGPG}).} Maximal reductive representations of $\Gamma$ to $SU(p,q)$, $p\leq q$, can be conjugated into $S(U(p,p)\times U(q-p))$. Maximal reductive representations of $\Gamma$ to $SO^{*}(2n)$, $n$ odd, can be conjugated into $SO^{*}(2n-2)\times SO(2)$.
\end{thm}

In turn, Theorem \ref{BGPG} is a special case of a more general result.

\begin{defi}
Say a Hermitian symmetric space is of {\em tube type} if it can be realized as a domain in $\bc^n$ of the form $\br^n +iC$ where $C\subset\br^n$ is a proper open cone.
\end{defi}

\begin{ex}
Siegel's upper half spaces and Grassmannians with isometry groups $PO(2,q)$ are of tube type.

The Grassmannian $\mathcal{D}_{p,q}$, $p\leq q$, with isometry group $PU(p,q)$ is of tube type iff $p=q$.

The Grassmannian $\mathcal{G}_{n}$ with isometry group $SO^* (2n)$ is of tube type iff $n$ is even.

The exceptional Hermitian symmetric space of dimension 27 is of tube type, the other one (of dimension 16) is not.

Products of tube type spaces are of tube type, so polydisks are of tube type.
\end{ex}

\begin{rem}
All maximal tube type subsymmetric spaces in a Hermitian symmetric space are conjugate. For instance, the maximal tube type subsymmetric space in $\mathcal{D}_{p,q}$ is $\mathcal{D}_{p,p}$. The maximal tube type subsymmetric space in $\mathcal{G}_{2n+1}$ is $\mathcal{G}_{2n}$.
\end{rem}

\begin{thm}
\label{BIW03}
{\em (Burger, Iozzi, Wienhard, \cite{BIW}).} Let $\Gamma$ be a closed surface group and $X$ a Hermitian symmetric space. Every maximal representation $\Gamma\to Isom(X)$ stabilizes a maximal tube type subsymmetric space $Y$. Conversely, for every tube type Hermitian symmetric space $X$, $Isom(X)$ admits Zariski dense maximal surface subgroups.
\end{thm}

\subsection{Results}

Our main result is a converse of Theorem \ref{BGPG} (i.e. Theorem \ref{BIW03} for classical simple Lie groups).

\begin{sthm}
\label{classical}
Let $G$ be a classical real Lie group, i.e. a real form of $SL(n,\bc)$, $O(n,\bc)$ or $Sp(n,\bc)$. Let $\Gamma$ be the fundamental group of a closed surface of genus $\geq 2\mathrm{dim}(G)^2$. A homomorphism $\phi:\Gamma\to G$ can be approximated by Zariski dense representations, unless the symmetric space $X$ of $G$ is Hermitian and not of tube type, and $\phi$ is maximal.
\end{sthm}
In other words, the exceptions are $G=SU(p,q)$, $q>p$ and
$\phi(\Gamma)$ is contained in a conjugate of $S(U(p,p)\times
U(q-p))\subset SU(p,q)$, or $G=SO^* (2n)$, $n$ odd, and
$\phi(\Gamma)$ is contained in a conjugate of $SO^* (2n-2)\times
SO(2)\subset SO^* (2n)$.

The genus assumption is probably unnecessary, but we are unable to remove it.

\begin{qu}
Does Theorem \ref{classical} extend to exceptional simple Lie groups ?
\end{qu}
Since flexibility easily holds in compact and complex groups (see Proposition \ref{cc}), and the rank one example $F_{4}^{-20}$ is treated in \cite{KP}, there remains 10 cases ($E_6^6$, $E_{6}^{-14}$, $E_{6}^{-26}$, $E_7^7$, $E_7^{-5}$, $E_7^{-25}$, $E_8^{8}$, $E_8^{-24}$, $F_4^4$, $G_2^2$) with at least one rigidity case ($E_6^{-14}$).

\subsection{Scheme of proof}

The proof relies on
\begin{itemize}
  \item a necessary and sufficient condition for flexibility from \cite{KP}, see Theorem \ref{flexssimple} below;
  \item tools from Burger, Iozzi and Wienhard's theory of tight maps between Hermitian symmetric spaces, \cite{BIW_tight};
  \item a detailed analysis of centers of centralizers of reductive subgroups of classical simple Lie groups.
\end{itemize}
This last analysis is performed in a rather bare handed manner, based on bilinear and sesquilinear algebra. This is where exceptional simple Lie groups elude us.

\subsection{Plan of the paper}

Section \ref{criteria} recalls the needed result from \cite{KP}. Section
\ref{tight} proves relevant consequences of the theory of tight
maps. Section \ref{list} provides a description of classical real simple
Lie groups as fixed points of involutions which helps in computing root
space decompositions in the complexified Lie algebra. This is done in
section \ref{zz} for $\mathfrak{sl}(n,\bc)$ and in section \ref{osu} for
$\mathfrak{so}(n,\bc)$ and $\mathfrak{sp}(n,\bc)$. The method consists in
first computing the root space decomposition in the standard representation
of the complexified Lie algebra, and deducing the decomposition in the
adjoint representation. The theory of tight homomorphisms allows to exclude
balancedness (with exceptions), first for real forms of $SL(n,\bc)$ in
section \ref{sl}, then for real forms of $SO(n,\bc)$ and $Sp(n,\bc)$ in
section \ref{so}. Theorem \ref{classical} is proven in section
\ref{nonreductive}.

\subsection{Acknowledgements}

Many thanks to Marc Burger and Jean-Louis Clerc, who explained us tight
maps and bounded symmetric domains.

\section{Flexibility criterion}
\label{criteria}

As far as the flexibility of a homomorphism $\phi$ is concerned, a key role is played by the center of the centralizer of the image of $\phi$. It splits the complexified Lie algebra of $G$ into root spaces $\mathfrak{g}_{\lambda}$. When the root $\lambda$ is pure imaginary, $\mathfrak{g}_{\lambda}$ carries a natural nondegenerate sesquilinear form defined as follows. Let $(X,X')\mapsto X\cdot X'$ denote the Killing form on $\mathfrak{g}\otimes\bc$. Then the sesquilinear form
\begin{eqnarray*}
s_{\lambda}(X,X')=\bar{X}\cdot X'
\end{eqnarray*}
 is nondegenerate on $\mathfrak{g}_{\lambda}$. Let $\Omega_{\lambda}$ denote the imaginary part of $s_{\lambda}$. It is a symplectic form on $\mathfrak{g}_{\lambda}$ viewed as a real vectorspace. The representation of $\Gamma$ on $\mathfrak{g}_{\lambda}$ gives rise to a homomorphism $\Gamma\to Sp(\mathfrak{g}_{\lambda},\Omega_{\lambda})$, an isometric action on the Siegel domain, and thus a Toledo invariant $T_{\lambda}$.

\begin{defi}
\label{defc}
Let $\mathfrak{t}\subset\mathfrak{g}$ be a torus, centralized by a homomorphism $\phi:\Gamma\to G$. Among the roots of the adjoint action of $\mathfrak{t}$ on $\mathfrak{g}$, let $P$ be the subset of pure imaginary roots $\lambda$ such that $2T_{\lambda}=-\chi(\Gamma)\mathrm{dim}_{\bc}(g_{\lambda})$, i.e. the symplectic $\Gamma$ action on $g_{\lambda}$ is maximal with positive Toledo invariant. Say $\mathfrak{t}$ is \emph{balanced} with respect to $\phi$ if $0$ belongs to the interior of the sum of the convex hull of the imaginary parts of elements of $P$ and the linear span of the real and imaginary parts of roots not in $\pm P$.
\end{defi}

Here is a necessary and sufficient condition for flexibility, for surface groups of sufficiently large genus.

\begin{thm}
\label{flexssimple}
{\em (\cite{KP}, Theorem 3)}.
Let $G$ be a semisimple real algebraic group. Let $\Gamma$ be the fundamental group of a closed surface of genus $\geq 2\mathrm{dim}(G)^2$. Let $\phi:\Gamma\to G$ be a homomorphism with reductive Zariski closure. Then $\phi$ is flexible if and only if $\mathfrak{c}$, the center of the centralizer of $\phi(\Gamma)$, is balanced with respect to $\phi$.
\end{thm}

The proof of Theorem \ref{classical} will rely on this criterion: we shall describe centers of centralizers of reductive subgroups and their roots and pile up restrictions on the set $P$ that make non balancedness exceptional. The reduction from nonreductive to reductive homomorphisms will be explained in section \ref{nonreductive}.

\section{Tightness}
\label{tight}

We collect in this section properties related to maximality of representations. We start with elementary facts. Deeper results will follow from tightness theory.

\subsection{Preservation of maximality}

We shall need that certain embeddings of Lie groups preserve maximality.

\begin{defi}
\label{defmaxpres} Let $\rho:H\to G$ be a homomorphism between
reductive Hermitian groups and $F:Y\to X$ denote an equivariant
totally geodesic map between the corresponding Hermitian symmetric
spaces. Let $\omega_X$ (resp. $\omega_Y$) denote the K\"ahler form,
normalized so that the minimum sectional curvature equals $-1$. Say
$\rho$ (or $F$) is {\em positively maximality preserving} if
\begin{eqnarray*}
\frac{1}{\mathrm{rank}(X)}F^{*}\omega_{X}=\frac{1}{\mathrm{rank}(Y)}\omega_{Y}.
\end{eqnarray*}
Say $\rho$ (or $F$) is merely {\em maximality preserving} if above equality holds
up to sign.
\end{defi}
Clearly, if this is the case, a homomorphism $\phi:\Gamma\to H$ is maximal
if and only if $\rho\circ\phi$ is. Furthermore, positively maximality
preserving maps do not change the signs of Toledo invariants.

\begin{lemma}
\label{max}
Isometric and holomorphic embeddings between equal rank Hermitian symmetric spaces are positively maximality preserving.
\end{lemma}

\begin{ex}
\label{max2}
The embedding $\mathcal{D}_{p,p}\hookrightarrow\mathcal{D}_{p,q}$ between Grassmannians, corresponding to the embedding $SU(p,p)\hookrightarrow SU(p,q)$, is isometric and holomorphic, and thus positively maximality preserving.
\end{ex}

\begin{ex}
\label{max3}
The embedding of symmetric spaces $\mathcal{G}_{n}\hookrightarrow\mathcal{G}_{n+1}$, corresponding to the embedding $SO^{*}(2n)\hookrightarrow SO^{*}(2n+2)$, is isometric and holomorphic, and thus, if $n$ is even, positively maximality preserving.
\end{ex}

\begin{ex}
\label{max4}
The embedding of symmetric spaces $\mathcal{S}_{n}\hookrightarrow\mathcal{S}_{n+1}$, corresponding to the embedding $Sp(2n,\br)\hookrightarrow Sp(2n+2,\br)$, is isometric and holomorphic, and thus never maximality preserving.
\end{ex}

We also need to understand when maximality is preserved under linear algebraic operations. Let us start with an easy case.

\begin{lemma}
\label{projective} Let $\Gamma$ be a surface group. Let $I_0$ be a
unitary $1$-dimensional representation of $\Gamma$ and $W$ a
sesquilinear representation of $\Gamma$. Then $Hom(I_0,W)$ is a
maximal representation if and only if $W$ is.
\end{lemma}

\begin{pf}
$Hom(I_0 ,W)$ and $W$ are isomorphic as projective representations. Since $U(W)$ acts on the symmetric space of $SU(W)$ via its quotient $PU(W)$, maximality is a projectively invariant property.
\end{pf}

But we shall need a more general case in subsection \ref{flexsu}. The following Lemma is a preparation for Lemma \ref{compuni}.

\begin{lemma}
\label{conjsum}
Let $\Gamma$ be a surface group. Let $W$, $W'$ be sesquilinear representations of $\Gamma$. Then
\begin{enumerate}
  \item $T(\bar{W})=-T(W)$.
  \item $T(W\oplus W')=T(W)+T(W')$.
\end{enumerate}
\end{lemma}

\begin{pf}
1. Passing from $W$ to $\bar{W}$ changes the sign of the complex structure on the symmetric space $X$ of $SU(W)$. This changes the sign of the K\"ahler form, and thus the sign of Toledo invariants.

2. Let $Y, Y'$ be the symmetric spaces of $SU(W)$ and $SU(W')$
respectively. When $Y\times Y'$ is mapped to the symmetric space $X$
of $SU(W\oplus W')$, the K\"ahler form $\omega_X$ of $X$ restricts
on the complex totally geodesic manifold $Y\times Y'$ to $\omega_Y
+\omega_{Y'}$, so Toledo invariants add up.
\end{pf}

\begin{lemma}
\label{compuni}
Let $\Gamma$ be a surface group. Let $V$ be a unitary representation of $\Gamma$ and $W$ a sesquilinear representation of $\Gamma$. Then $Hom(V,W)$ is a maximal representation of $\Gamma$ if and only if $W$ is.
\end{lemma}

\begin{pf}
Let $\phi:\Gamma\to U(Hom(V,W))$ denote the Hom of the two given representations. Changing the given action on $V$ into the trivial representation gives rise to a representation $\psi:\Gamma\to U(Hom(V,W))$. Let us show that $\phi$ and $\psi$ have equal Toledo invariants. Split $W=W^+ \oplus W^-$ into positive definite and negative definite subspaces. Then $Hom(V,W^+)$ is a maximal positive definite subspace of $Hom(V,W)$, i.e. a point in the Hermitian symmetric space $X$ associated to $SU(Hom(V,W))$. Its $U(V)\times U(W)$-orbit $Y$ is totally geodesic in $X$, and thus contractible. Therefore, one can choose an equivariant map $\tilde{f}:\tilde{\Sigma}\to X$ whose image is contained in $Y$. Since $U(V)$ fixes $Y$ pointwise, $\tilde{f}$ is equivariant with respect to both $\phi$ and $\psi$. Therefore the corresponding Toledo invariants are the same. Since $\psi$ is a direct sum of $\mathrm{dim}(V)$ copies of the action on $W$, Lemma \ref{conjsum} applies, so $T(Hom(V,W))=\mathrm{dim}(V)T(W)$. Since $\mathrm{rank}(X)$ equals $\mathrm{dim}(V)$ times the rank of the symmetric space associated to $SU(W)$, $\psi$, and thus $\phi$, is maximal if and only if the original $\Gamma$ action on $W$ is.
\end{pf}

\subsection{Tightness}

Tightness theory is a way to draw strong consequences from the existence of maximal representations.

\begin{defi}
\label{deftight}
{\em (Burger, Iozzi, Wienhard, \cite{BIW_tight})}.
Let $G$ be a {\em reductive Hermitian} group, i.e. a connected reductive Lie group in which the center is compact and such that the symmetric spaces $X_i$ associated to all simple noncompact factors is Hermitian. Normalize the metric on $X_i$ so that the minimum holomorphic sectional curvature equals $-1$. Let $\kappa_{G}^{b}$ denote the bounded continuous cohomology class on $G$ defined by integrating the K\"ahler form of $X=\prod X_i$ on triangles with geodesic sides. Let $\Gamma$ be locally compact group. Say a continuous homomorphism $\phi:\Gamma\to G$ is {\em tight} if
\begin{eqnarray*}
\parallel \phi^* \kappa_{G}^{b}\parallel =\parallel \kappa_{G}^{b}\parallel.
\end{eqnarray*}
\end{defi}

\begin{ex}
Maximal homomorphisms of surface groups to reductive Hermitian groups are tight.
\end{ex}

\subsection{Maximality preserving versus tight}

Here is the basic mechanism which makes tightness enter our arguments: if $\phi:\Gamma\to G$ is maximal and factors through $\rho:H\to G$, then $\rho$ is tight. There is a converse statement.

\begin{prop}
\label{maxprestight}
Let $H$, $G$ be reductive Hermitian groups. Assume that the symmetric space associated to $H$ is irreducible. Let $\rho:H\to G$ be a continuous homomorphism. Then $\rho$ is maximality preserving if and only if $\rho$ is tight.
\end{prop}

\begin{pf}
This follows from Proposition 2.12 of \cite{BIW_tight}.
\end{pf}

\begin{ex}
\label{exmax}
(Example 8.7 of \cite{BIW_tight}). The obvious embeddings $SU(n,n)\to Sp(4n,\br)$ and $SO^{*}(4n)\to Sp(8n,\br)$ are tight. It follows that $SO^{*}(4n)\to SU(2n,2n)$ is tight. All three embeddings are thus maximality preserving.
\end{ex}

Direct proofs of these facts will be given in an appendix, Lemmas
\ref{max1} and \ref{sosu}.

\subsection{Consequences of tightness}

\begin{lemma}
\label{tighttight}
\begin{enumerate}
  \item Let $H\subset G$ be connected real algebraic groups. If $G$ is reductive Hermitian and the embedding $H\hookrightarrow G$ is tight, then $H$ is reductive Hermitian too.
  \item Let $\rho:G\to G'$ be a tight homomorphism between reductive Hermitian groups. If the kernel of $\rho$ is compact and $G'$ is of tube type, so is $G$.
\end{enumerate}
\end{lemma}

\begin{pf}
This is a combination of Theorems 7.1 and 6.2 of \cite{BIW_tight}.
\end{pf}

\begin{lemma}
\label{lemsupq}
Let $V$, $V'$ be vectorspaces equipped with nondegenerate sesqui\-linear forms. Assume that $Hom(V,V')$, equipped with the natural sesquilinear form
\begin{eqnarray*}
(f,f')\mapsto\mathrm{Trace}(f^{*}\circ f'),
\end{eqnarray*}
(here, $f^{*}$ denotes the adjoint with respect to the sesquilinear forms on $V$ and $V'$), has vanishing signature. Assume that the induced homomorphism $U(V)\times U(V')\to U(Hom(V,V'))$ is tight (see Definition \ref{deftight}). Then one of $V$ and $V'$ is definite and the other has vanishing signature.
\end{lemma}

\begin{pf}
Since $U(Hom(V,V'))$ is Hermitian of tube type, Lemma \ref{tighttight} implies that $U(V)\times U(V')$ is of tube type, up to compact groups. This implies that each of the sesquilinear vectorspaces $V$ and $V'$ is either of vanishing signature or definite. Clearly, if both are definite, $U(V)\times U(V')$ is compact and tightness is out of sight.

Let us show that $V'$ and $V'$ cannot both have vanishing signature, i.e. one of them must be definite. For this, we use the tightness criterion of \cite{BIW_tight}, Corollary 8.2. Let $J_{\ell}$ denote the generator of the center of the maximal compact subgroup of $U(V)$ which defines the complex structure on the symmetric space of $U(V)$. In a splitting of $V=V^{+}\oplus V^{-}$ in a sum of positive (resp. negative) definite subspaces,
\begin{eqnarray*}
J_{V}=\begin{pmatrix}
\frac{i}{2}I_{\frac{d_{V}}{2}} & 0\\
0  & -\frac{i}{2}I_{\frac{d_{V}}{2}}
\end{pmatrix}
\end{eqnarray*}
Use a similar splitting $V'=V'^{+}\oplus V'^{-}$ and split accordingly $Hom(V,V')$, whose elements become block matrices $N=\begin{pmatrix}
A & B\\
C & D
\end{pmatrix}$. The action of $U_{V}\times U_{V'}$ is $(U,U',N)\mapsto U'^{-1}NU$. The linearized action is $(U,U',N)\mapsto -U'N+NU$. The image of $(J_{V},J_{V'})$ under the linearized action is
$N\mapsto \begin{pmatrix}
0 & -iB   \\
iC  & 0
\end{pmatrix}$, a map whose matrix is $\begin{pmatrix}
0 & 0 & 0 & 0\\
0 & 0 & 0 & 0\\
0 & 0 & -i & 0\\
0 & 0 & 0 & i
\end{pmatrix}$
in a splitting of $Hom(V,V')$ into a positive subspace $Hom(V^{+},V'^{+})\oplus Hom(V^{-},V'^{-})$ and a negative subspace $Hom(V^{+},V'^{-})\oplus Hom(V^{-},V'^{+})$. Therefore the inner product with the element
\begin{eqnarray*}
J_{Hom(V,V')}=\begin{pmatrix}
\frac{i}{2}I_{\frac{d_{V}}{2}}  & 0 & 0 & 0\\
0 & \frac{i}{2}I_{\frac{d_{V}}{2}}  & 0 & 0\\
0 & 0 & -\frac{i}{2}I_{\frac{d_{V}}{2}}  & 0\\
0 & 0 & 0 & -\frac{i}{2}I_{\frac{d_{V}}{2}}
\end{pmatrix}
\end{eqnarray*}
vanishes. According to \cite{BIW_tight}, Corollary 8.2, the homomorphism
$U(V)\times U(V')\ra U(Hom(V,V'))$ is not tight. We conclude that one of $V$ and $V'$ has vanishing signature and the other one is definite.
\end{pf}

The following example of non tight embedding will be useful in the proof of Corollary \ref{fin}.

\begin{lemma}
\label{o22}
The injection $\iota:O(2,2)\hookrightarrow U(2,2)$ is not tight.
\end{lemma}

\begin{pf}
$O(2,2)\subset U(2,2)$ is the fixed point set of conjugation $\sigma$, i.e. $\sigma\circ\iota=\sigma$. $\sigma$ induces an orientation reversing isometry of the symmetric space $X$ of $U(2,2)$, which changes the sign of the K\"ahler form, $\sigma^{*}\omega_{X}=-\omega_{X}$. $\iota$ induces a totally geodesic embedding of symmetric spaces, still denoted by $\iota:Y\to X$. Since $\iota^{*}\omega_{X}=\iota^{*}\sigma^{*}\omega_{X}=-\iota^{*}\omega_{X}$, $\iota^{*}\omega_{X}=0$, $\iota^{*}$ is not isometric in bounded continuous cohomology, so $\iota$ is not tight.
\end{pf}

\subsection{Alternative definition of balancedness}

One may replace symplectic structures by sesquilinear structures in the definition of balancedness.

\begin{prop}
\label{maxtight}
Let $G$ be a semisimple real algebraic group. Let
$\Gamma$ be a surface group.  Let $\phi:\Gamma\to G$ be a
homomorphism. Let $\mathfrak{c}$ denote the center of the
centralizer of $\phi(\Gamma)$. Let $Z_{G}(\mathfrak{c})$ denote its centralizer in $G$. Let $\lambda$ be a pure imaginary
root of the adjoint action of $\mathfrak{c}$ on
$\mathfrak{g}\otimes \bc$. Let $\mathfrak{g}_{\lambda}$ denote
the corresponding root space, equipped with a sesquilinear form $s_{\lambda}$ and its imaginary part, the symplectic form $\Omega_{\lambda}$. Then the symplectic representation $\Gamma\to Sp(g_{\lambda},\Omega_{\lambda})$ is maximal with positive Toledo invariant if and only if
\begin{enumerate}
  \item $s_{\lambda}$ has vanishing signature;
  \item the sesquilinear representation $\Gamma\to U(g_{\lambda},s_{\lambda})$ is maximal with positive Toledo invariant.
  \end{enumerate}
If this is the case, then the homomorphism $Z_{G}(\mathfrak{c})\to U(g_{\lambda},s_{\lambda})$ is tight.
\end{prop}

\begin{pf}
1. Assume that $\Gamma\to Sp(g_{\lambda},\Omega_{\lambda})$ is maximal. It factors through $U(g_{\lambda},s_{\lambda})$. Thus the inclusion $U(g_{\lambda},s_{\lambda})\hookrightarrow Sp(g_{\lambda},\Omega_{\lambda})$ is tight. Lemma \ref{tighttight} implies that $U(g_{\lambda},s_{\lambda})$ has tube type, i.e. the signature of $s_{\lambda}$ vanishes.

2. According to Example \ref{exmax}, the embedding $U(g_{\lambda},s_{\lambda})\hookrightarrow Sp(g_{\lambda},\Omega_{\lambda})$ is positively maximality preserving, so maximality and positivity of Toledo invariant do not change when passing from symplectic to unitary groups.

3. The maximal representation $\Gamma\to U(g_{\lambda},s_{\lambda})$ factors via the homomorphism $Z_{G}(\mathfrak{c})\to U(g_{\lambda},s_{\lambda})$, which must be tight itself.
\end{pf}

Now we can explain how balancedness will be analyzed in the sequel. For classical simple Lie groups, root spaces $\mathfrak{g}_{\lambda}$ turn
out to be expressible as mapping spaces $Hom(V,V')$, centralizers
$Z_{G}(\mathfrak{c})$ are products of classical simple Lie
groups. Proposition \ref{maxtight}, combined with Lemmas \ref{tighttight} and \ref{lemsupq}, restricts the possible groups involved, as will be seen in sections \ref{sl} and \ref{so}.

\section{Classical simple Lie groups}
\label{list}

\subsection{Definition}
\label{defcla}

Classical simple real Lie groups are special linear groups of division rings $D$ with center $\br$ and special unitary groups of nondegenerate binary forms over $D$, see the appendix in \cite{S}.

There are only 3 such division rings: $\br$, $\bc$ and $\bh$, leading to 3 special linear groups, $SL(n,\br)$, $SL(n,\bc)$, and $SL(n,\bh)$.

Let $\iota$ be a continuous (anti)-automorphism of $D$: $\iota$ can
be identity, complex conjugation or quaternionic conjugation. Let
$\epsilon=\pm 1$. By a $(\iota,\epsilon)$-symmetric
binary\footnote{To avoid confusion, we keep the words Hermitian for
positive definite forms, and sesquilinear for the complex case.}
form on a right $D$-vectorspace $V$, we mean a $D$-valued
$\br$-bilinear map $h:V\times V\to D$ such that
\begin{enumerate}
  \item for all $v$, $v'\in V$ and all $q\in D$, $h(v,v'q)=h(v,v')q$;
  \item for all $v$, $v'\in V$, $h(v',v)=\epsilon\iota(h(v,v'))$.
\end{enumerate}
Note that from (2), it follows that $H(vq,v')=\bar q H(v,v')$. Say
$h$ is {\em nondegenerate} if the only $v\in V$ such that
$h(v,v')=0$ for all $v'\in V$ is $0$.

The group $U(V,h)$ consists of right $D$-linear self maps of $V$ which preserve $h$. We are interested in the special unitary group $SU(V,h)=U(V,h)\cap SL(V,D)$.

When $D=\br$, the only choice for $\iota$ is identity, leading to
\begin{itemize}
  \item if $\epsilon=1$, real special orthogonal groups, indexed by dimension and signature and denoted by $SO(p,q)$, where $p+q=\mathrm{dim}_{\br}(V)$;
  \item if $\epsilon=-1$, real symplectic groups, indexed by dimension only and denoted by $Sp(n,\br)$, where $n=\mathrm{dim}_{\br}(V)$ is even.
\end{itemize}

When $D=\bc$ there are two choice for $\iota$. If $\iota$ is identity, this leads to complex special orthogonal (if $\epsilon=1$) and symplectic (if $\epsilon=-1$) groups, indexed by dimension and denoted respectively by $SO(n,\bc)$ and $Sp(n,\bc)$. If $\iota$ is complex conjugation, $\epsilon=\pm1$ lead to groups called special unitary groups, indexed by dimension and signature and denoted by $SU(p,q)$, where $p+q=\mathrm{dim}_{\bc}(V)$.

When $D=\bh$ there are two choices for $\iota$. If $\iota$ is identity, no nonzero $(\iota,\epsilon)$-symmetric forms exist. If $\iota$ is quaternionic conjugation, this leads to
\begin{itemize}
  \item if $\epsilon=1$, quaternionic unitary groups, indexed by dimension and signature and denoted by $Sp(p,q)$, $p+q=\mathrm{dim}_{\bh}(V)$;
  \item if $\epsilon=-1$, quaternionic skew-unitary groups, indexed by dimension only and denoted by $SO^* (2n)$, where $n=\mathrm{dim}_{\bh}(V)$ is even.
\end{itemize}

\subsection{Real forms of complex groups}
\label{rf}

We shall be mainly concerned with the 7 families of non complex groups, 3 attached to binary forms on real vectorspaces, 1 on complex vectorspaces and 3 on quaternionic vectorspaces. Each of these groups is obtained as the fixed point set of an anti-$\bc$-linear involutive automorphism $\sigma$ of a complex Lie group, as we now explain.

Consider first the complex unitary family $U(p,q)$. Here, the data is a nondegenerate sesquilinear form $s$ on a complex vectorspace $V$. For $f\in GL(V)$, let $f^{*}$ denote the $s$-adjoint of $f$, defined by
\begin{eqnarray*}
\forall v,\,v'\in V,\quad s(f(v),v')=s(v,f^{*}(v')).
\end{eqnarray*}
Then $\sigma(f)=(f^{*})^{-1}$ is an anti-$\bc$-linear involutive automorphism of $GL(V)$ and of $SL(V)$. The fixed point set of $\sigma$ in $SL(V)$ is $SU(V,s)$.

\medskip

The 6 remaining families admit a common construction.

Given a real vectorspace $V_{\br}$, let $V=V_{\br}\otimes\bc$ and let $\tau=$conjugation. Note that $\tau^{-1} =\tau$. If $b$ is a nondegenerate quadratic or symplectic form on $V_{\br}$, let $B=b\otimes\bc$.

Given a right quaternionic vectorspace $V_{\bh}$, pick a basis $(1,i,j,k)$ of $\bh$, use right multiplication by $i$ to turn $V_{\bh}$ into a complex vectorspace denoted by $V$. Let $\tau$ be right multiplication by $j$. Note that $\tau^{-1} =-\tau$. For $q=a+jb\in\bh$, $a$, $b\in\bc$, denote by $\mathcal{C}(q)=a$. If $h$ is a nondegenerate $(\bar{\hskip1em},\epsilon)$-symmetric binary form on $V_{\bh}$, let, for $v$, $v'\in V$,
\begin{eqnarray*}
B(v,v')=\mathcal{C}(h(vj,v')).
\end{eqnarray*}
Then $B$ is a nondegenerate (-$\epsilon$)-symmetric $\bc$-bilinear form on $V$.

In both cases, $\tau$ is anti-$\bc$-linear, $\tau^{-1}=\eta\tau$ for some $\eta\in\pm1$, and
\begin{eqnarray*}
B(\tau(v),\tau(v'))=\overline{B(v,v')}.
\end{eqnarray*}
For $f\in End_{\bc}(V)$, set
\begin{eqnarray*}
\sigma(f)=\tau\circ f\circ\tau^{-1}.
\end{eqnarray*}
Then $\sigma$ is an anti-$\bc$-linear involutive automorphism of the
algebra $\mathfrak{gl}(V)=End_{\bc}(V)$ and of its subgroups $SL(V)$
and $O(V,B)$. The fixed point set of $\sigma$ in $SL(V)$ is
$SL(V_{\br})$ (resp. $SL(V_{\bh})$). The fixed point set of $\sigma$
in $O(V,B)$, i.e., $\tau\circ f=f\circ \tau$, is $U(V_{\br},b)$
(resp. $U(V_{\bh},h)$). This construction yields 6 of the 7 families
of non complex classical groups.

\subsection{Consequences for roots}

The special form of the involution $\sigma$ for 6 of the 7 families of noncomplex groups has the following consequence.

\begin{lemma}
\label{conjugate}
Let $G$ belong to one of the above 6 families of classical simple Lie groups. Let $H\subset G$ be a reductive subgroup. Let $\mathfrak{c}\subset\mathfrak{g}$ be the center of its centralizer. Let $\ell$ be a root of $\mathfrak{c}$ on $\bc^n$. Let $I_{\ell}$ denote the corresponding root space. Then $\bar{\ell}$ is a root and $I_{\bar{\ell}}=\tau(I_{\ell})$. Furthermore,
\begin{enumerate}
  \item if $G=SL(n,\br)$, $O(p,q)$ or $Sp(n,\br)$, $I_{\bar{\ell}}=\overline{I_{\ell}}$;
  \item if $G=SL(n,\bh)$, $SO^* (2n)$ or $Sp(p,q)$, $I_{\ell}+ I_{\bar{\ell}}$ is a quaternionic subspace;
  \item If $\ell'\notin\{\ell,\bar{\ell}\}$, $I_{\ell}+ I_{\bar{\ell}}$ and $I_{\ell'}+ I_{\bar{\ell'}}$ are orthogonal with respect to the $(\iota,\epsilon)$-symmetric binary form $h$, and therefore nondegenerate.
\end{enumerate}
\end{lemma}

\begin{pf}
Let $Z\in\mathfrak{c}$ and $v\in I_{\ell}$. Then
$$
Z(\tau(v))=\tau(Z(v))=\tau(\ell(Z)v)=\overline{\ell(Z)}\tau(v),
$$
showing that $I_{\bar{\ell}}=\tau(I_{\ell})$. For the 3 real families, $\tau$ is conjugation, thus $I_{\bar{\ell}}=\overline{I_{\ell}}$. For the 3 quaternionic families, $\tau$ is right multiplication by $j$, so $I_{\ell}+ I_{\bar{\ell}}$ is stable by right multiplication by $i$ and $j$, i.e. a quaternionic subspace.

Since $h$ is $G$-invariant, for all $v$, $v'\in V$ and $Z\in\mathfrak{c}$, $h(Z(v),v')+h(v,Z(v'))=0$. If $v\in I_{\ell}$ and $v'\in I_{\ell'}$,
\begin{eqnarray*}
0&=&h(Z(v),v')+h(v,Z(v'))\\
&=&h(v\ell(Z),v')+h(v,v'\ell'(Z))\\
&=&h(v,v')\overline{\ell(Z)}+h(v,v')\ell'(Z)\\
&=&h(v,v')(\bar{\ell}+\ell')(Z).
\end{eqnarray*}
This implies that $h(v,v')=0$ if $\ell'\notin\{\ell,\bar{\ell}\}$.
\end{pf}

\subsection{Killing form}

Here is a formula for the sesquilinear structure appearing in Theorem \ref{flexssimple}, valid in all cases.

\begin{lemma}
\label{killing} Let $\mathfrak{g}^{\bc}\subset \mathfrak{sl}(n,\bc)$
be a complex Lie subalgebra of $\mathfrak{sl}(n,\bc)$. Let $\sigma$
be an anti-$\bc$-linear involutive automorphism of
$\mathfrak{g}^{\bc}$, with fixed point set $\mathfrak{g}$. The
natural sesquilinear form $s$ on $\mathfrak{g}\otimes\bc$ arising
from the Killing form of $\mathfrak{g}$ (see section \ref{criteria})
is proportional to
\begin{eqnarray*}
s(X,X')=\mathrm{Trace}(\sigma(X)\circ X').
\end{eqnarray*}
\end{lemma}

\begin{pf}
Since $\sigma$ is anti-$\bc$-linear, for $f\in\mathfrak{g}$, $\sigma(if)=-i\sigma(f)=-if$. It follows that the $-1$-eigenspace of $\sigma$ in $\mathfrak{g}^{\bc}$ is $i\mathfrak{g}$. The map
\begin{eqnarray*}
\mathfrak{g}^{\bc}\to\mathfrak{g}\otimes\bc,\quad X\mapsto(\frac{X+\sigma(X)}{2},\frac{X-\sigma(X)}{2i})
\end{eqnarray*}
is an isomorphism. It pulls back conjugation on $\mathfrak{g}\otimes\bc$ to $\sigma$ on $\mathfrak{g}^{\bc}$, thus it pulls back the sesquilinear form $\bar{X}\cdot X$ of $\mathfrak{g}\otimes\bc$ to $\sigma(X)\cdot X$ on $\mathfrak{g}^{\bc}$. The Killing form on $\mathfrak{g}^{\bc}$ is proportional to the restriction to $\mathfrak{g}^{\bc}$ of the Killing form of $\mathfrak{sl}(n,\bc)$, whence the formula $\mathrm{Trace}(\sigma(X)\circ X)$.
\end{pf}

\subsection{Flexibility of compact and complex Lie groups}

For completeness' sake, we recall here the treatment of compact and
complex Lie groups from \cite{KP}.

\begin{prop}
\label{cc}
Let $G$ be a compact or complex semisimple Lie group. Let $\mathfrak{c}\subset\mathfrak{g}$ be the center of the centralizer of a reductive subgroup. Then $\mathfrak{c}$ is balanced.
\end{prop}

\begin{pf}
Let $X\in\mathfrak{g}\otimes\bc$, $X=f+ig$. Then $s(X,X)=\bar{X}\cdot X=f\cdot f+g\cdot g$. If $G$ is compact, the Killing form is negative definite, so is $s$. For all roots $\lambda$ of $\mathfrak{c}$ in the adjoint representation, the sesquilinear form $s_{\lambda}$ on the root space $\mathfrak{g}_{\lambda}$ is negative definite. No root has vanishing signature. According to Proposition \ref{maxtight}, $P$ is empty so $\mathfrak{c}$ is balanced.

If $G$ is complex, centralizers are complex Lie subgroups, $\mathfrak{c}$ is a complex vectorsubspace, and roots $\lambda$ are $\bc$-linear maps. None of them is pure imaginary (i.e. takes pure imaginary values on $\mathfrak{c}$). Thus $P$ is empty so $\mathfrak{c}$ is balanced.
\end{pf}

\section{Complexified centers of centralizers in \texorpdfstring{$SL(n,\bc)$}{}}
\label{zz}

The first step is to list the possible complexified centers $\mathfrak{c}\otimes\bc$ and describe the root structure. This depends only on the complexified Lie algebra.

\begin{lemma}
\label{centrcompl}
Let $\mathfrak{g}$ be a real Lie algebra. Let $H\subset\mathfrak{g}$ be a subset. Let $ZZ_{\mathfrak{g}}(H)$ denote the center of its centralizer. Then
\begin{eqnarray*}
ZZ_{\mathfrak{g}}(H)\otimes\bc=ZZ_{\mathfrak{g}\otimes\bc}(H).
\end{eqnarray*}
\end{lemma}

Our approach consists in using the standard complex representation $V$ of $\mathfrak{g}\otimes\bc$. The root space decomposition of $V$ under $\mathfrak{c}\otimes\bc=ZZ_{\mathfrak{g}\otimes\bc}(H)$ is obtained from the isotypical decomposition of $V$ under $H$. This easily provides us with the decomposition of $End(V)$ under $\mathfrak{c}\otimes\bc$, and then of the invariant subspace $\mathfrak{g}\otimes\bc\subset End(V)$.

In this section, we treat the model case of $SL(n,\bc)$, and in the next section, the more elaborate cases $O(n,\bc)$ and $Sp(n,\bc)$.

\subsection{\texorpdfstring{$H$}{}-modules}

\begin{defi}
Let $H$ be a group. The data of a finite dimensional complex vectorspace and a homomorphism of $H$ onto a reductive real algebraic subgroup of $Gl(V)$ is called an {\em $H$-module}.
\end{defi}

Note that $H$-invariant subspaces of $H$-modules are again $H$-modules. A $H$-module is {\em irreducible} if it has no proper $H$-submodules.

\begin{defi}
\label{defirr}
Let $E(H)$ denote the set of equivalence classes of $H$-modules. Given an $H$-module $V$ and $\pi\in E(H)$, let $I_{\pi}$, the $\pi$-isotypical component of $V$, denote the sum of all submodules of $V$ belonging to the equivalence class $\pi$.
\end{defi}

The following Lemma is well known, but we give a full proof since it serves as a model for orthogonal, symplectic and unitary versions of it given in the next section.

\begin{lemma}
\label{0}
Any $H$-module $V$ splits as a direct sum of its isotypical components
\begin{eqnarray*}
V=\bigoplus_{\pi\in E(H)}I_{\pi}.
\end{eqnarray*}
\end{lemma}

\begin{pf}
Let $W\subset I_{\pi}$ be an irreducible invariant subspace. Since $H$ is reductive, for each invariant subspace $Z$ of $I_{\pi}$ belonging to the class $\pi$, there exists an $H$-invariant complement to $Z$, thus an $H$-invariant projector $p_{Z}:I_{\pi}\to Z$. Since such submodules $Z$ generate $I_{\pi}$, for at least one such $Z$, $p_{Z}(W)\not=0$, thus $p_{Z}(W)=Z$ and $W$ belongs to $\pi$.

More generally, if $F\subset E(H)$ is a subset, and $W\subset \sum_{\pi\in F}I_{\pi}$ is an irreducible invariant subspace, then $W$ belongs to one of the classes in $F$. Indeed, otherwise $p_{Z}(W)=0$ for all invariant subspaces $Z$ whose class belongs to $F$, and these generate $\sum_{\pi\in F}I_{\pi}$. In particular, for every $\pi\in E(H)$, $I_{\pi}\cap\sum_{\pi'\not=\pi}I_{\pi'}=\{0\}$, which shows that the sum $\sum_{\pi\in E(H)}I_{\pi}$ is direct.

Since $H$ is reductive, the invariant subspace $\bigoplus_{\pi\in E(H)}I_{\pi}$ admits an invariant complement, which contains an irreducible subspace. This is a contradiction unless $\bigoplus_{\pi\in E(H)}I_{\pi}=V$.
\end{pf}

\subsection{Centers of centralizers in \texorpdfstring{$SL(n,\bc)$}{}}
\label{zzsln}

\begin{lemma}
\label{centrslnr}
Let $H\subset SL(n,\bc)$ be a reductive subgroup. Let $\mathfrak{c}\subset\mathfrak{sl}(n,\bc)$ be the center of its centralizer. Let $L$ denote the set of nonzero roots of $\mathfrak{c}$ in the standard representation of $\mathfrak{sl}(n,\bc)$, and $d_\ell$ the dimensions of the corresponding eigenspaces. Then $L$ has $\mathrm{dim}(\mathfrak{c})+1$ elements, which satisfy exactly one linear relation,
$$\sum_{\ell\in L}d_\ell \ell=0.$$
Furthermore, the map
$$L\times L\setminus diagonal\ra \Lambda,\quad (\ell,\ell')\mapsto \ell-\ell'$$
is one to one onto the set $\Lambda$ of nonzero roots of $\mathfrak{c}$ in the adjoint representation of $\mathfrak{sl}(n,\bc)$.
\end{lemma}

\begin{pf}
Since $H$ is reductive, the standard action of $H$ on $\bc^n$ splits into irreducibles. Let us group them into isotypical components $I_{\ell}$: $I_{\ell}$ is the direct sum of $k_{\ell}$ isomorphic irreducible summands. According to Schur's Lemma, the group of $H$-automorphisms of the representation $I_{\ell}$ is isomorphic to $GL(k_{\ell},\bc)$. Then
\begin{eqnarray*}
Z_{GL(n,\bc)}(H)&=&\prod_{\ell\in L}Z_{GL(I_{\ell})}(H_{|I_{\ell}})\\
&\simeq&\prod_{\ell\in L}GL(k_{\ell},\bc),
\end{eqnarray*}
whose center is $(\bc^{*})^{L}$, acting on $\bc^n$ by multiplication
by a different constant on each $I_{\ell}$. Pick a basis of $\bc^n$
adapted to the splitting $\bc^n =\bigoplus_{\ell}I_{\ell}$. Then the
center $\mathfrak{c}$ of $Z_{SL(n,\bc)}(H)$ consists of diagonal
matrices $diag(a_{1},\ldots,a_{n})$ whose entries corresponding to
basis vectors from the same $I_{\ell}$ are equal, and which sum up
to $0$. It follows that the elements of $L$ generate $\mathfrak c^*$
and satisfy only one linear relation, $\sum_{\ell\in
L}d_{\ell}\ell=0$. In particular, if $(\ell,\ell')$, $(m,m')$ are
distinct ordered pairs of distinct elements of $L$,
$\ell-\ell'-m+m'$ does not vanish identically on $\mathfrak{c}$.
This shows that the map
$$L\times L\setminus diagonal\ra \Lambda,\quad (\ell,\ell')\mapsto \ell-\ell',$$
which is clearly surjective onto the set of nonzero roots of the adjoint action on $\mathfrak{sl}(n,\bc)$, is injective as well.
Furthermore, the root space for $\ell-\ell'$ is $Hom(I_{\ell},I_{\ell'})$.
\end{pf}

\section{Centralizers in orthogonal, symplectic or unitary groups}
\label{osu}

In this section, the rootspace decomposition of $\mathfrak{g}$ under the center of the centralizer of a reductive subgroup is given, when $\mathfrak{g}=\mathfrak{so}(n,\bc)$ or $\mathfrak{sp}(n,\bc)$. This is a first step in handling real forms of these Lie algebras. With little extra effort, one can treat simultaneously the case of $\mathfrak{g}=\mathfrak{su}(p,q)$. This will help treating this particular real form of $\mathfrak{sl}(n,\bc)$.

Let $(V,B)$ be a complex vectorspace equipped with a nondegenerate binary form $B$ of one of the following three types,
\begin{enumerate}
  \item symmetric bilinear,
  \item skew-symmetric bilinear,
  \item symmetric sesquilinear,
\end{enumerate}
which we denote by $(\iota,\epsilon)$-symmetric, $\epsilon=\pm 1$, $\iota=1$ or $\bar{\hskip1em}$ (identity or conjugation). Note that the combination $(\iota,\epsilon)=(\bar{\hskip1em},-1)$ makes perfect sense but does not bring anything new, since if $B$ is a skew-symmetric sesquilinear form, then $iB$ is symmetric sesquilinear. Let $O=O^{\iota,\epsilon}(V,B)$ denote its automorphism group (note that if $B$ is skew-symmetric bilinear (resp. sesquilinear), this is a symplectic (resp. unitary) rather than an orthogonal group, whence the notation $O^{\iota,\epsilon}$). Let $H\subset O$ be a reductive real algebraic subgroup. In this section, we describe the center of the centralizer of $H$ in $O$.

\subsection{Bilinear and sesquilinear forms}

\begin{nota}
Let $\iota$ be a continuous automorphism of $\bc$, i.e. either conjugation or identity. Let $\epsilon\in\{-1,1\}$. Let $V$ be a complex vectorspace. A {\em $(\iota,\epsilon)$-symmetric form} on $V$ is a real bilinear form $B:V\times V\to \bc$ such that
\begin{itemize}
  \item for $\lambda$, $\lambda'\in\bc$, $v$, $v'\in V$, $B(\lambda v,\lambda' v')=\iota(\lambda)\lambda' B(v,v')$;
  \item $B(v',v)=\epsilon \iota B(v,v')$.
\end{itemize}
\end{nota}
In other words, if $\iota=1$ and $\epsilon=1$, $B$ is symmetric bilinear. If $\iota=1$ and $\epsilon=-1$, $B$ is skew-symmetric bilinear. If $\iota=\bar{\hskip1em}$ and $\epsilon=1$, $B$ is symmetric sesquilinear. If $\iota=\bar{\hskip1em}$ and $\epsilon=-1$, $B$ is skew-symmetric sesquilinear. We shall ignore the fourth case, since if $B$ is a skew-symmetric sesquilinear form, then $iB$ is symmetric sesquilinear.

\begin{nota}
1. Let $\iota$ be a continuous automorphism of $\bc$, i.e. either conjugation or identity. Let $V$ be a complex vectorspace. Then $V^{\iota}$ means $V$ if $\iota=1$, $\bar{V}$ if $\iota=\bar{\hskip1em}$. Also, $V^{\iota,*}$ means the dual vectorspace $V^*$ if $\iota=1$, $\bar{V}^*$ (i.e. the space of anti-$\bc$-linear forms on $V$) if $\iota=\bar{\hskip1em}$.

2. Let $B$ be a $(\iota,\epsilon)$-symmetric form on $V$. Let $\sharp_B :V\to V^{\iota,*}$ denote the $\bc$-linear map which maps $v\in V$ to the anti-$\bc$-linear functional
\begin{eqnarray*}
\sharp_{B}(v):v'\mapsto B(v',v).
\end{eqnarray*}
Say that $B$ is {\em nondegenerate} if $\sharp_B$ is an isomorphism. If not, its kernel is called the {\em nullspace} of $B$.
\end{nota}

\begin{nota}
Given a linear map $L:V\to V^{\iota,*}$, there is an adjoint map $L^{\iota,\top}:V\to V^{\iota,*}$ defined by
\begin{eqnarray*}
\langle L^{\iota,\top}(v'),v \rangle=\langle L(v),v'\rangle,
\end{eqnarray*}
where $\langle v^{*},v \rangle$ denotes the evaluation of a linear or anti-linear form $v^{*}$ on a vector $v\in V$.
\end{nota}

If $L$ is $\epsilon$-symmetric, i.e. if
$L^{\iota,\top}=\epsilon\iota L$, the formula
\begin{eqnarray*}
B(v,v')= \langle L(v'),v\rangle
\end{eqnarray*}
defines a $(\iota,\epsilon)$-symmetric form such that $\sharp_B =L$. Therefore $B$ and $\sharp_B$ are equivalent data.

\begin{ex}
\label{exsesq}
Let $W$ be a complex vectorspace. Then the tautological isomorphism of $V=W\times W^{\iota,*}$ to $V^{\iota,*}$ gives rise to a tautological $(\iota,\epsilon)$-symmetric form on $V$,
\begin{eqnarray*}
(v,v^*)\cdot(w,w^*)= \epsilon\iota(\langle v^* ,w\rangle)+\langle w^* ,v\rangle.
\end{eqnarray*}
\end{ex}

Every two nondegenerate $(1,\epsilon)$-symmetric forms are isomorphic. On the other hand, $(\iota,1)$-symmetric forms, i.e. symmetric sesquilinear forms, take real values on the diagonal, so sign and signature issues arise: two nondegenerate $(\iota,1)$-symmetric forms are isomorphic if and only if they have the same signature. For instance, the tautological form of Example \ref{exsesq} has vanishing signature.

\subsection{Bilinear and sesquilinear \texorpdfstring{$H$}{}-modules}

\begin{defi}
Let $H$ be a group. The data of a finite dimensional complex vectorspace equipped with a bilinear (either symmetric or skew-symmetric) or sesquilinear form and a homomorphism of $H$ onto a reductive real algebraic subgroup of its automorphism group will be called a {\em $(\iota,\epsilon)$-linear $H$-module}.
\end{defi}

\begin{lemma}
\label{lem1}
Let $W$ be an irreducible $H$-module. The space of $H$-invariant bilinear (resp. sesquilinear) forms on $W$ has dimension at most $1$. A non-zero $H$-invariant bilinear (resp. sesquilinear) form is automatically non-degenerate, and in the bilinear case, it is either symmetric or skew-symmetric.
\end{lemma}

\begin{pf}
Let $b$ be an $H$-invariant bilinear form on $W$. Its nullspace is $H$-invariant. Therefore $b$ is either zero or non-degenerate. Assume $b$ is nonzero and denote by $\sharp_b : W\to W^{\iota,*}$ the corresponding isomorphism (in the sesquilinear case, $\iota=$ conjugation, $W^{\iota,*}=\bar{W}^*$). Let $b'$ be an other $H$-invariant bilinear form on $W$. Then $L=(\sharp_{b})^{-1}\circ\sharp_{b'}$ is an $H$-equivariant endomorphism of $W$, thus $L$ is a multiple of identity (Schur's Lemma). This shows that $b'$ is a multiple of $b$. $b$ has a symmetric and a skew-symmetric component. They have to be linearly dependant, in the bilinear case ($\iota=1$), this implies that one of them vanishes. Therefore $b$ is either symmetric or skew-symmetric.
\end{pf}

\begin{defi}
\label{deforth}
Say an $H$-module is {\em $(\iota,\epsilon)$-orthogonal} if it admits an invariant non-degenerate $(\iota,\epsilon)$-symmetric form.
\end{defi}

\begin{cor}
\label{class} The classification of irreducible
$(\iota,\epsilon)$-linear $H$-modules can be deduced from the
classification of irreducible $H$-modules: the forgetful map
$E^{(\iota,\epsilon)}(H)\to E(H)$ is onto, the fiber of an
equivalence class of irreducible $H$-modules contains $1$ or $2$
elements depending wether it is $(\iota,\epsilon)$-orthogonal or
not.
\end{cor}

\begin{lemma}
\label{2} Let $(V,B)$ be a $(\iota,\epsilon)$-linear $H$-module. Let
$W$ and $W'$ be distinct irreducible $H$-invariant subspaces. Assume
that $W'$ is not orthogonal to $W$. Then $W'$ is isomorphic, as an
$H$-module, to the (conjugate-)dual $W^{\iota,*}$ of $W$ (in the
sesquilinear case, $W^{\iota,*}=\bar{W}^*$).
\end{lemma}

\begin{pf}
The map $v\mapsto (\sharp_B (v))_{|W}$, $W'\to W^{\iota,*}$, is $H$-equivariant. According to Schur's Lemma, such a map is either zero or an isomorphism, and all such maps are proportional. By assumption, it does not vanish, thus $W'$ and $W^{\iota,*}$ are isomorphic $H$-modules.
\end{pf}

\begin{defi}
Let $H$ be a group, $V$ an $H$-module. Say $V$ is {\em bi-isotypical} if their exists an irreducible $H$-module $Z$ such that every irreducible
invariant subspace $W\subset V$ is isomorphic either to $Z$ or to $Z^{\iota,*}$.
\end{defi}

\begin{cor}
\label{3}
Let $(V,B)$ be a non-degenerate $(\iota,\epsilon)$-linear $H$-module. Then $V$ canonically splits as an orthogonal direct sum of its bi-isotypical components,
\begin{eqnarray*}
V=\bigoplus_{\pi}I_{\pi,\pi^{\iota,*}},
\end{eqnarray*}
where, given an equivalence class $\pi$ of irreducible $H$-modules, $I_{\pi,\pi^{\iota,*}}=I_{\pi}+I_{\pi^{\iota,*}}$ is the sum of all irreducible invariant subspaces of $V$ isomorphic either to $\pi$ or to $\pi^{\iota,*}$. Furthermore, $I_{\pi,\pi^{\iota,*}}$ is non-degenerate, and the centralizer $Z_{O(V,B)}(H)$ of $H$ in the automorphism group is isomorphic to a direct product,
\begin{eqnarray*}
Z_{O(V,B)}(H)=\prod_{\pi}Z_{O(I_{\pi,\pi^{\iota,*}})}(H).
\end{eqnarray*}
\end{cor}

\subsection{Examples of bi-isotypical bilinear/sesquilinear \texorpdfstring{$H$}{}-modules}

From now on, we analyze non-degenerate bi-isotypical $(\iota,\epsilon)$-linear $H$-modules. Here are two examples.

\begin{ex}
\label{exdual}
Let $Z=\pi^{\oplus r}$ be an isotypical $H$-module. Set $V=Z\times Z^{\iota,*}$, equip it with the canonical $(\iota,\epsilon)$-symmetric form
\begin{eqnarray*}
(v,v^*)\cdot(w,w^*)= \epsilon\iota(\langle v^* ,w\rangle)+\langle w^* ,v\rangle.
\end{eqnarray*}
Then the centralizer of $H$ in $O^{(\iota,\epsilon)}(V,\cdot)$ is
isomorphic to an orthogonal, symplectic or unitary group, if $\pi$
and $\pi^{\iota,*}$ are equivalent
\begin{eqnarray*}
Z_{O^{(\iota,\epsilon)}(V,\cdot)}(H)&\cong& O^{(\iota,\epsilon)}(2r,\bc)\\
&:=&\begin{cases}
O(2r,\bc) & \text{ if }\iota=1,~\epsilon=1, \\
Sp(2r,\bc) & \text{ if }\iota=1,~\epsilon=-1, \\
U(r,r)& \text{ if }\iota=\bar{\hskip1em},~\epsilon=1,
\end{cases}
\end{eqnarray*}
otherwise to a general linear group,
\begin{eqnarray*}
Z_{O^{(\iota,\epsilon)}(V,\cdot)}(H)\cong Gl(r,\bc).
\end{eqnarray*}
\end{ex}

\begin{pf}
Pick a basis $e_i$ of $W$ and take $r$ copies of it to form a basis
of $Z$.  Take the image of this basis under $\sharp_{B}$ to get a
basis of $Z^{\iota,*}$, i.e.  for $\epsilon=1$, choose
$B(e_i,e_i)=1$ and $e_i^*=\sharp_{B}(e_i)$, and then take
$\{e_1,\cdots; e_1^*,\cdots\}$ as a basis for $V$. For
$\epsilon=-1$, choose $e_1, \cdots, e_{2k}$ so that
$B(e_i,e_{i+k})=1=-B(e_{i+k},e_i)$ and $B(e_i,e_j)=0$ otherwise.
Choose a basis for $V$ in this case as
$\{e_1,\cdots,e_k,e_{k+1},\cdots,e_{2k},e_{k+1}^*,\cdots,e_{2k}^*,e_1^*,\cdots,e_k^*\}$.
This gives a basis of $V$ in which the matrix of the
bilinear/sesquilinear form $B$ equals $\begin{pmatrix}
0 & 1 \\
\epsilon  & 0
\end{pmatrix}$ (blocks have size $rd$ where $d=\mathrm{dim}\,\pi$). In this basis, the matrix of an element $g$ of $H$ splits into blocks of size $d$, with the first $r$ diagonal blocks equal to $\pi(g)$ and the last $r$ equal to $\pi(g^{-1})^{\iota,\top}$, all other blocks vanish.

If $\pi$ and $\pi^{\iota,*}$ are equivalent, endomorphisms of $V$ which commute with $H$ have matrices whose blocks of size $d$ are scalar, i.e. proportional to the unit $d\times d$ matrix. In other words, they can be written $A\otimes 1$ where $A\in Gl(2r,\bc)$. Such a matrix preserves $B$ if and only if
\begin{eqnarray*}
(A\otimes 1)^{\iota,\top}(\begin{pmatrix}
0 & 1 \\
\epsilon  & 0
\end{pmatrix}\otimes 1)(A\otimes 1)=\begin{pmatrix}
0 & 1 \\
\epsilon & 0
\end{pmatrix}\otimes 1,
\end{eqnarray*}
i.e. iff $A$ belongs to $O(2r,\bc)$ if $\epsilon=1$, to $Sp(2r,\bc)$ if $\epsilon=-1$ and to $U(r,r)$ if $\iota=$conjugation. In other words,
\begin{eqnarray*}
Z_{O^{(\iota,\epsilon)}(V,\cdot)}(H)\cong O^{(\iota,\epsilon)}(2r,\bc).
\end{eqnarray*}

If $\pi$ and $\pi^{\iota,*}$ are not equivalent, endomorphisms of $V$ which commute with $H$ preserve the splitting $V=I_{\pi}\oplus I_{\pi^{\iota,*}}$ and have matrices whose blocks of size $d$ are scalar, i.e. proportional to the unit $d\times d$ matrix. In other words, they can be written $(A\oplus A')\otimes 1$ where $A$, $A'\in Gl(r,\bc)$. Such a matrix preserves $B$ if and only if $A'=(A^{\iota,\top})^{-1}$. In other words,
\begin{eqnarray*}
Z_{O^{(\iota,\epsilon)}(V,\cdot)}(H)\cong Gl(r,\bc).
\end{eqnarray*}
\end{pf}

\begin{ex}
\label{exsum}
Let $V=\pi^{\oplus r}$ be an isotypical $H$-module such that $\pi$ is $(\iota,\epsilon)$-orthogonal, i.e. preserves a non-degenerate $(\iota,\epsilon)$-symmetric form $b$. Let $D$ be a real diagonal invertible $r\times r$ matrix. Set $B=b\otimes D$, i.e. $(V,B)$ is an orthogonal direct sum of $r$ real multiples of the same non-degenerate $(\iota,\epsilon)$-linear $H$-module. Then the centralizer of $H$ in $O^{\epsilon}(V,b\otimes D)$ is isomorphic to an orthogonal/unitary group,
\begin{eqnarray*}
Z_{O^{\iota,\epsilon}(V,b\otimes D)}(H)\cong O^{\iota}(\bc^{r},D).
\end{eqnarray*}
\end{ex}

Note that in the bilinear case, $O(\bc^{r},D)=O(r,\bc)$ is a genuine orthogonal group, even when $\epsilon=-1$, i.e. when we deal with skew-symmetric forms. In the sesquilinear case, $O^{\iota}(\bc^{r},D)=U(p,q)$ where $p-q=\mathrm{sign}(D)$.

\begin{pf}
Repeat the same basis of $\pi$ to get a basis of $V$. The matrix of $B=b\otimes D$ is diagonal in blocks of size $d=\mathrm{dim}\,\pi$, with diagonal blocks equal to real multiples of $b$. Element $g\in H$ acts by a diagonal matrix in $d\times d$-blocks, with diagonal blocks equal to $\pi(g)$. The centralizer of $H$ in $Gl(V)$ consists of matrices with scalar $d\times d$-blocks, i.e. of the form $A\otimes 1$ for $A\in Gl(r,\bc)$. $(\iota,\epsilon)$-orthogonal matrices satisfy
\begin{eqnarray*}
(A\otimes 1)^{\iota,\top}(b\otimes D)(A\otimes 1)=b\otimes(A^{\iota,\top}DA)=b\otimes D,
\end{eqnarray*}
i.e. $A^{\iota,\top}DA=D$. Thus
\begin{eqnarray*}
Z_{O^{\epsilon}(V,B)}(H)\cong O^{\iota}(\bc^{r},D).
\end{eqnarray*}
\end{pf}

\subsection{Classification of bi-isotypical bi/sesquilinear \texorpdfstring{$H$}{}-modules}

There are 3 cases, depending wether $\pi$ is $(\iota,\epsilon)$-orthogonal, $(\iota,-\epsilon)$-orthogonal, or neither. In each case, we will need the following lemma.

\begin{lemma}
\label{isot}
Let $V=W\oplus W'$ be a $(\iota,\epsilon)$-linear $H$-module where $W$, $W'$ are isotropic, isotypic, and the map $L:W'\to W^{\iota,*}$, $L(v')=(\sharp_{B} (v'))_{|W}$ is an isomorphism. Then $V$ is isomorphic to Example \ref{exdual}.
\end{lemma}

\begin{pf}
The $H$-map $(w+w')\mapsto (w,L(w'))$ is an isometry $W\oplus W' \to W\times W^{\iota,*}$, since
\begin{eqnarray*}
(v,L(v'))\cdot(w,L(w'))
&=&\epsilon\iota(\langle L(v'),w\rangle) +\langle L(w'),v\rangle \\
&=&\epsilon\iota(\langle \sharp_{B} (v'),w\rangle) +\langle \sharp_{B} (w'),v\rangle\\
&=&\epsilon\iota(B(w,v'))+B(v,w')\\
&=&B(v',w)+ B(v,w')\\
&=&B(v+v',w+w').
\end{eqnarray*}
\end{pf}

\begin{prop}
\label{non}
Let $(V,B)$ be a non-degenerate bi-isotypical bilinear $H$-module. Assume $\pi$ and $\pi^{\iota,*}$ are not isomorphic. Then $(V,B)$ is isomorphic to Example \ref{exdual}.
\end{prop}

\begin{pf}
Irreducible invariant subspaces of $V$ belong either to $\pi$ or $\pi^{\iota,*}$, which do not admit non-degenerate invariant bilinear/sesquilinear forms, thus all are isotropic. Lemma \ref{2} implies that every two distinct irreducible invariant subspaces of $I_{\pi}$ are orthogonal. Thus $I_{\pi}$ is isotropic, and so is $I_{\pi^{\iota,*}}$. Since $V$ is non-degenerate, $\sharp_{B}:V\to V^{\iota,*}$ induces isomorphisms $L:I_{\pi}\to (I_{\pi^{\iota,*}})^{\iota,*}$ and $L':I_{\pi^{\iota,*}}\to(I_{\pi})^{\iota,*}$ related by $L^{\bot}=\epsilon L'$. According to Lemma \ref{isot}, this shows that $V$ is isomorphic, as a $(\iota,\epsilon)$-linear $H$-module, to Example \ref{exdual}.
\end{pf}

\begin{prop}
\label{-eps}
Let $V=I_{\pi}$ be an isotypical $H$-module equipped with a nondegenerate $\epsilon$-symmetric bilinear form $B$. Assume $\pi$ is $(-\epsilon)$-orthogonal. Then $(V,B)$ is isomorphic to Example \ref{exdual}.
\end{prop}

\begin{pf}
Since $\pi$ is not $\epsilon$-orthogonal, all irreducible $H$-submodules of $V$ are isotropic. Let $W$ be one of them. A $H$-invariant complement to $W^{\bot}$ contains an irreducible $H$-submodule $W'$, which is not orthogonal to $W$. The $H$-map $L:W'\to W^*$, $v'\mapsto (\sharp_{B} (v'))_{|W}$ is non-zero, thus an isomorphism, and Lemma \ref{isot} implies that $W\oplus W'$ is isomorphic to Example \ref{exdual}, in particular, it is non-degenerate. Its orthogonal is again non-degenerate, isotypic, modelled on a $(-\epsilon)$-orthogonal $H$-module. By induction on dimension, $V$ is an orthogonal direct sum of copies of Example \ref{exdual}, thus isomorphic to Example \ref{exdual}.
\end{pf}

\begin{prop}
\label{pm}
Let $V=I_{\pi}$ be an isotypical $H$-module equipped with a nondegenerate $(\iota,\epsilon)$-symmetric bilinear form $B$. Assume $\pi$ is $(\iota,\epsilon)$-orthogonal. Then $(V,B)$ is isomorphic to Example \ref{exsum}.
\end{prop}

\begin{pf}
Let us show that $V$ contains at least one non-degenerate irreducible invariant subspace. Pick an irreducible $H$-submodule $W$ of $V$. If it is non-degenerate, we are done. Otherwise, $W$ is isotropic. An $H$-invariant complement to $W^{\bot}$ contains an irreducible $H$-submodule $W'$, which is not orthogonal to $W$. If $W'$ is non-degenerate, we are done. Otherwise, $W'$ is isotropic too. Then $W\oplus W'$ is isomorphic to $W\times W^{\iota,*}$ equipped with the canonical $(\iota,\epsilon)$-symmetric form. Indeed, the $H$-map $L:W'\to W^{\iota,*}$, $v'\mapsto (\sharp_{B} (v'))_{|W}$ is non-zero, thus an isomorphism, and Lemma \ref{isot} applies. By assumption, there exists an $\epsilon$-symmetric $H$-isomorphism $M:W\to W^{\iota,*}$. Then the graph $Z=\{(w,M(w))\,|\,w\in W\}$ of $M$ is non-degenerate. Indeed, for $v$, $w\in W$,
\begin{eqnarray*}
(v,M(v))\cdot(w,M(w))&=&\epsilon\iota(\langle M(v),w\rangle) +\langle M(w),v\rangle\\
&=&2\epsilon (\langle M(v),w\rangle)
\end{eqnarray*}
cannot vanish for all $w$, unless $v=0$. $Z$ is a
non-degenerate irreducible $(\iota,\epsilon)$-linear $H$-submodule of $W\times W^{\iota,*}$, which embeds
isometrically into $V$, this is the required subspace.

The proof of the proposition is concluded by induction on dimension. If $V$ is irreducible, we are done. Otherwise, we just showed that $V$ has at least one non-degenerate irreducible submodule, say $W$. Then the induction hypothesis applies to its orthogonal $W^{\bot}$.
\end{pf}

\begin{cor}
\label{centralizorth}
Let $(V,B)$ be a nondegenerate bi-isotypical $(\iota,\epsilon)$-linear $H$-module. Then $(V,B)$ is isomorphic either to Example \ref{exdual} or to Example \ref{exsum}. The centralizer of $H$ in the automorphism group of $(V,B)$ is isomorphic to a general linear group in the former case, an orthogonal group in the latter.
\end{cor}

\subsection{Centers of centralizers in unitary groups}

In this case, we directly get information on a real form.

\begin{prop}
\label{centrcentrunit}
Let $(V,B)$ be a non-degenerate symmetric sesquilinear $H$-module. Let $\mathfrak{c}$ denote the center of the centralizer of $H$ in $U(V,B)$. Then root spaces of $V$ under $\mathfrak{c}$ correspond to isotypical components $I_{\pi}$ under $H$. They fall into bi-isotypical components $BI_{\pi}$ which are pairwise orthogonal. If $I_{\pi}=BI_{\pi}$, then the corresponding root of $\mathfrak{c}$ is pure imaginary. If $I_{\pi}\not=BI_{\pi}$, then $BI_{\pi}=I_{\pi}\oplus I_{\bar{\pi}^{*}}$, the corresponding roots $\ell_{\pi}$ and $\ell_{\pi^{*}}$ are opposite. The linear relations between roots are generated by $\ell_{\pi}+\ell_{\bar{\pi}^{*}}=0$, $\pi\not=\bar{\pi}^{*}$.

Let $N$ be a set which contains exactly one element of each pair $\{\pi,\bar{\pi}^{*})\}$ of equivalence classes of irreducible sesquilinear $H$-modules occurring in $V$, such that $\pi\not=\bar{\pi}^{*}$. Then the map $(\ell_{\pi})_{\pi\in L_N}:\mathfrak{c}\to\bc^{L_N}$ is onto.
\end{prop}

\begin{pf}
If $I_{\pi}=BI_{\pi}$, then, as a sesquilinear $H$-module, $BI_{\pi}$ is isomorphic to Example \ref{exsum} and contributes a unitary factor to the centralizer of $H$, whose center is a pure imaginary subgroup of $SL(n,\bc)$. It acts on $BI_{\pi}$ by multiplication by a pure imaginary number. Thus $BI_{\pi}$ is a root space for a pure imaginary root. Otherwise, $BI_{\pi}$ is isomorphic to Example \ref{exdual} and contributes a general linear group factor to the centralizer. Its center is a complex subgroup of $SL(n,\bc)$, it acts on $I_{\pi}$ (resp. $I_{\bar{\pi}^{*}}$) by multiplication by an unrestricted complex number (resp. the opposite number). This produces a subspace $\mathfrak{c}'$ of $\mathfrak{c}$ which admits a complex structure, the corresponding roots are $\bc$-linear and half of them (to avoid the relations $\ell_{\pi}+\ell_{\bar{\pi}^{*}}=0$) provide complex coordinates on $\mathfrak{c}'$.
\end{pf}

\begin{cor}
\label{rootssu}
Let $\tilde{L}$ be the set of roots, $L_I$ the subset of roots which take only pure imaginary values. Then for every $\ell \in L\setminus L_I$, $-\ell$ is again a root. Let $L_N$ be a set which contains exactly one element of each pair $\{\ell,-\ell\}$, $\ell\in L\setminus L_I$. Then
\begin{itemize}
  \item $I_{\ell}$ and $I_{\ell'}$ are orthogonal unless $\ell'=-\ell$;
  \item if $\ell\in L_I$, $I_{\ell}$ is nondegenerate;
  \item the map $(\ell)_{\ell\in L_N}:\mathfrak{c}\to\bc^{L_N}$ is onto.
\end{itemize}
\end{cor}

\subsection{Centers of centralizers in \texorpdfstring{$\epsilon$}{}-orthogonal groups}

We continue our convention that
\begin{eqnarray*}
O^{\epsilon}(n,\bc)=\begin{cases}
O(n,\bc) & \text{ if }\epsilon=1, \\
Sp(n,\bc) & \text{ if }\epsilon=-1.
\end{cases}
\end{eqnarray*}
Note that in the latter case, $n$ has to be even.

\begin{prop}
\label{centrcentrorth}
Let $(V,B)$ be a non-degenerate bilinear $H$-module. Let $\mathfrak{c}$ denote the center of the centralizer of $H$ in $O^{\epsilon}(V,B)$. Under $\mathfrak{c}$, $V$ splits into root spaces as follows. Each isotypical component $I_{\pi}$ where $\pi$ is an irreducible $H$-module which is not equivalent to its contragredient $\pi^*$ is a root space for a non-zero root $\ell_{\pi}$. It is isotropic. The sum of all other isotypical components constitutes the $0$ root space. Bi-isotypical components $I_{\pi}+I_{\pi^{*}}$ are pairwise orthogonal. Relations among non-zero roots are generated by the following
\begin{eqnarray*}
\ell_{\pi}+\ell_{\pi^{*}}=0.
\end{eqnarray*}
Thus the number of non zero roots is $2\mathrm{dim}\,\mathfrak{c}$. Let $L$ be a set containing exactly one element of each pair $(\pi,\pi^*)$. Then $(\ell_{\pi})_{\pi\in L}:\mathfrak{c}\to\bc^{L}$ is a linear bijection.
\end{prop}

\begin{pf}
Let $L'$ be the set of equivalence classes of irreducible $H$-modules which are isomorphic to their contragredient. Then
\begin{eqnarray*}
Z_{O^{\epsilon}(V,B)}(H)=\prod_{\pi\in
L}GL(r_{\pi},\bc)\times\prod_{\pi\in L'}O^{\epsilon}(2r_{\pi},\bc),
\end{eqnarray*}
thus
\begin{eqnarray*}
ZZ_{O^{\epsilon}(V,B)}(H)=\bigoplus_{\pi\in L}\bc(id_{I_{\pi}}-id_{I_{\pi^{*}}}),
\end{eqnarray*}
i.e. $(\ell_{\pi})_{\pi\in L}:\mathfrak{c}\to\bc^{L}$ is a linear bijection. Furthermore, if $\pi\in L$, the bi-isotypical component $I_{\pi}\oplus I_{\pi^{*}}$ splits into two root spaces relative to roots $\ell_{\pi}$ and $\ell_{\pi^{*}}=-\ell_{\pi}$.
\end{pf}

\begin{cor}
\label{rootorth}
Let $\tilde{L}$ denote the set of roots of $\mathfrak{c}$ on $\bc^n$. The roots of $\mathfrak{c}$ in its adjoint action on $\mathfrak{so}^{\epsilon}(V,B)$ are exactly all differences $\ell-\ell'$, for $\ell$, $\ell'\in \tilde{L}$, including $2\ell$ if $\epsilon=-1$ or $\mathrm{dim}(I_{\ell})>1$, excluding $2\ell$ if $I_{\ell}$ is $1$-dimensional and $\epsilon=1$.

In other words, if $L$ is a set of representatives of pairs $\{-\ell,\ell\}$ of nonzero roots, one finds $0$, all sums $\pm\ell\pm\ell'$ for distinct $\ell$, $\ell'\in L$, sometimes $\pm 2\ell$ (depending on $\epsilon$ and $\mathrm{dim}(I_{\ell})$), and, if $0$ is also a root of $\mathfrak{c}$ on $\bc^n$, all $\pm\ell$, $\ell\in L$.
\end{cor}

\begin{pf}
The Lie algebra $\mathfrak{so}(V,B)$ is the space of {\em $B$-skew-symmetric} endomorphisms of $V$, i.e. $\bc$-linear maps $f:V\to V$ satisfying, for all $v$, $v'\in V$,
\begin{eqnarray*}
B(f(v),v')+B(v,f(v'))=0.
\end{eqnarray*}

The roots of $\mathfrak{c}$ in the adjoint representation are differences $\lambda=\ell-\ell'$ of roots of $\mathfrak{c}$ in $\bc^n$. If $\ell'\not=\pm\ell$, the root space relative to $\ell-\ell'$ is the subspace of $B$-skew-symmetric elements of $Hom(I_{\ell},I_{\ell'})\oplus Hom(I_{-\ell'},I_{-\ell})$. It does never vanish. Indeed, for every $f\in
Hom(I_{\ell},I_{\ell'})$, there is a unique $g\in Hom(I_{-\ell'},I_{-\ell})$ such that $(f,g)$ is $B$-skew-symmetric. The formula for $g$ is
\begin{eqnarray*}
g=-\sharp_{B}^{-1}\circ f^{\top}\circ\sharp_{B}.
\end{eqnarray*}

Here is an alternative description of $\mathfrak{so}(V,B)$: mapping
$f\in End(V)$ to the bilinear form $b(v,v')=B(f(v),v')$ identifies
$End(V)$ with the space $V^* \otimes V^*$ of $\bc$-bilinear forms on
$V$ and $\mathfrak{so}(V,B)$ with the subspace
$\Lambda^{\epsilon}V^{*}$ of (-$\epsilon$)-symmetric $\bc$-bilinear
forms on $V$. The adjoint action of $Z\in \mathfrak{so}(V)$ is
\begin{eqnarray*}
Zb(v,v')=B([Z,f](v),v')=-b(v,Zv')-b(Zv,v').
\end{eqnarray*}
If $v\in I_{\ell}$, $v'\in I_{\ell'}$ and $Z\in\mathfrak{c}$, then
\begin{eqnarray*}
Zb(v,v')=-(\ell+\ell')(Z)b(v,v').
\end{eqnarray*}
Therefore the root space relative to $2\ell$ identifies with $\Lambda^{\epsilon}I_{\ell}^{*}$. It vanishes if and only if $\epsilon=1$ and $\mathrm{dim}(I_{\ell})=1$.
\end{pf}

\section{Real forms of \texorpdfstring{$SL(n,\bc)$}{}}
\label{sl}

\subsection{Flexibility in \texorpdfstring{$SL(n,\br)$}{} and \texorpdfstring{$SL(n,\bh)$}{}}

\begin{prop}
\label{slnr}
Let $\mathfrak{c}\subset\mathfrak{sl}(n,\br)$ (resp. $\mathfrak{sl}(n/2,\bh)$) be the center of the centralizer of a reductive subgroup of $SL(n,\br)$ (resp. $SL(n/2,\bh)$). Let $\lambda$ be a pure imaginary root of $\mathfrak{c}$ in its adjoint action on $\mathfrak{sl}(n,\br)$. Then the signature of the Killing form restricted to $\mathfrak{g}_{\lambda,\br}$ does not vanish. It follows that $\mathfrak{c}$ is balanced.
\end{prop}

\begin{pf}
Let $\bc^n$ denote the standard representation of
$\mathfrak{sl}(n,\bc)$. As in subsection \ref{rf}, let
$\tau(v)=\bar{v}$ in the complex case, and $\tau(v)=vj$ in the
quaternionic case (here, $\bc^n =\bh^{n/2}$ is viewed as a right
quaternionic vectorspace). Let, for $f\in \mathfrak{sl}(n,\bc)$,
$\sigma(f)=\tau\circ f\circ\tau^{-1}$. Then
$Fix(\sigma)=\mathfrak{g}=\mathfrak{sl}(n,\br)$ (resp.
$\mathfrak{sl}(n/2,\bh)$).

Under $\mathfrak{c}$, $\bc^n$ splits into root spaces $\bc^n=\bigoplus_{\ell}I_{\ell}$, $\mathrm{dim}(I_{\ell})=d_{\ell}$. Roots are either real or come in pairs $\{\ell,\bar{\ell}\}$ (Lemma \ref{conjugate}).
According to Lemma \ref{centrslnr}, every nonzero root $\lambda$ of $\mathfrak{c}$ in its adjoint action on $\mathfrak{sl}(n,\bc)$ can be uniquely written in the form $\ell-\ell'$. Such a root is pure imaginary if and only if $\ell'=\bar{\ell}$, i.e. $\lambda=\ell-\bar{\ell}$. The corresponding root space is
$$
\mathfrak{g}_{\lambda}=Hom(I_{\ell},I_{\bar{\ell}}).
$$
Let $f\in Hom(I_{\ell},I_{\bar{\ell}})$, $\sigma(f)=\bar{f}\in Hom(I_{\bar{\ell}},I_{\ell})$. Pick a basis of $I_{\ell}$ and take its image by $\tau$ as a basis of $I_{\bar{\ell}}$. Let $M$ denote the matrix of $f$ in the chosen basis of $\bc^n$. Then the matrix of $\sigma(f)$ is $\sigma(M)$, and $\mathrm{Trace}(\sigma(f)\circ f)=\mathrm{Trace}(\bar{M}M)$. Write $M=S+A$ where $S$ is symmetric and $A$ is skew-symmetric. Then
\begin{eqnarray*}
\mathrm{Trace}(\bar{M}M)=\mathrm{Trace}(\bar{S}S)+\mathrm{Trace}(\bar{A}A)=\mathrm{Trace}(S^{*}S)-\mathrm{Trace}(A^{*}A)
\end{eqnarray*}
has signature $\mathrm{dim}(\{S\})-\mathrm{dim}(\{A\})=\frac{d_{\ell}(d_{\ell}+1)}{2}-\frac{d_{\ell}(d_{\ell}-1)}{2}=d_{\ell}$,
which is nonzero.
\end{pf}

\subsection{Flexibility in \texorpdfstring{$SU(p,q)$}{}}
\label{flexsu}

\begin{prop}
\label{supq} Let $\Gamma$ be a surface group, let $\phi:\Gamma\to
SU(p,q)$ a reductive homorphism, let $\mathfrak{c}$ be the center of
the centralizer of $\phi(\Gamma)$. Assume that $\mathfrak{c}$ is not
balanced with respect to $\phi$. Then, up to conjugacy,
$\phi(\Gamma)$ is contained in $S(U(p,p)\times U(q-p))$, and $\phi$
is maximal.
\end{prop}

\begin{pf}
Under $\mathfrak{c}$, the standard representation of
$\mathfrak{sl}(n,\bc)$ splits into root spaces,
$\bc^{n}=\bigoplus_{\ell}I_{\ell}$. The roots of $\mathfrak{c}$ in
the adjoint representation are differences $\ell-\ell'$, and all of
them indeed occur. Following Corollary \ref{rootssu}, split the set
of roots as $\tilde{L}=L_I \cup L_N \cup -L_N$, where $L_I$ is the
subset of pure imaginary roots. If $\ell\in L_N$, $2\ell$ is a root
in the adjoint representation, and it is not pure imaginary. If
$\ell\in L_I$ and $\ell'\in L_N$, $\ell-\ell'$ is not pure
imaginary. If $\ell$ and $\ell'\in L_N$, $\ell-\ell'$ is not pure
imaginary either, since it factors through a surjective map
$\mathfrak{c}\to\bc^{L_N}$ and a $\bc$-linear form
$\bc^{L_N}\to\bc$. So if $L_N$ is non empty, non pure imaginary
roots span $\mathfrak c^{*}$, and $\mathfrak{c}$ is balanced,
contradiction. So $L_N$ is empty, all roots on $\bc^n$ are pure
imaginary, and the corresponding root spaces $I_{\ell}$ are
nondegenerate and pairwise orthogonal. In the sequel, we shall
replace roots by there imaginary parts without expressing it in the
notation.

Each root $\lambda$ in the adjoint representation can be expressed in a unique way as $\lambda=\ell-\ell'$, and
$$
\mathfrak{g}_{\lambda}=Hom(I_{\ell},I_{\ell'}).
$$
A calculation shows that the signature of the natural Hermitian form on $Hom(I_{\ell},I_{\ell'})$
equals $-\mathrm{sign}(I_{\ell})\mathrm{sign}(I_{\ell'})$. This time, the signature is not automatically nonzero.
So different arguments, based on \cite{BIW_tight}, are needed.

Assume that the sesquilinear action of $\Gamma$ on
$\mathfrak{g}_{\lambda}$ is maximal. The situation we are
considering is as follows: since $H=\overline{\phi(\Gamma)}$
preserves $I_\ell$ and preserves a sesquilinear form on it,
$\phi(\Gamma)\subset U(I_\ell)$, and similarly for $I_{\ell'}$,
hence $\phi(\Gamma)\rightarrow U(I_\ell)\times U(I_{\ell'})$ and we
obtain
$$
\xymatrix{
\phi(\Gamma)\subset Z_G(\mathfrak c) \ar@{^(->}[r]^{}\ar@{^(->}[rd]_{f} & U(I_{\ell})\times U(I_{\ell'}) \ar[d]^{} \\
&U(\mathfrak g_\lambda)=U(Hom(I_\ell, I_{\ell'}))=U(I_\ell\otimes
I_{\ell'}^*)} \;\;,
$$ where $f$ is tight.
According to Proposition \ref{maxtight}, the sesquilinear space
$Hom(I_{\ell},I_{\ell'})$ must have vanishing signature. Also, the
homomorphism $Z_{G}(\mathfrak{c})\to U(Hom(I_{\ell},I_{\ell'}))$
must be tight.  Lemma \ref{lemsupq} implies that one of $I_{\ell}$
and $I_{\ell'}$ has vanishing signature and the other one is
definite. Say $I_{\ell}$ is definite, for instance. According to
Lemma \ref{compuni}, maximality of $Hom(I_{\ell},I_{\ell'})$ implies
maximality of the $\Gamma$ action on $I_{\ell'}$, with
$T_{\lambda}=\mathrm{dim}(I_{\ell})T(I_{\ell'})$. If instead
$I_{\ell'}$ is definite,
$T_{\lambda}=-T(I_{\ell})\mathrm{dim}(I_{\ell'})$.

Let $D$ (resp. $E$, resp. $O$) denote the set of roots $\ell$ such that $I_{\ell}$ is definite (resp. has vanishing signature, resp. has non vanishing signature). As in definition \ref{defc}, let $P$ denote the set of roots $\lambda$ such that $\rho_{\lambda}$ is maximal with positive Toledo invariant. Equivalently, of differences $\ell-\ell'$, $\ell\in D$, $\ell'\in E$, such that $I_{\ell'}$ is maximal with positive Toledo invariant. Let $N$ be the complement of $\pm P$.

Let us show that if $O$ is non empty, $\mathfrak{c}$ is balanced with respect to $\phi$. Indeed, let $\ell_0 \in O$. Then for all $\ell'\not=\ell_0$, for $\lambda=\ell_0 -\ell'$, $\rho_{\lambda}$ is not maximal, thus $\ell_0 -\ell'\in N$. Since the roots $\ell$ span $\mathfrak{c}^{*}\otimes\bc$ and satisfy the extra equation $\sum_{\ell}d_{\ell}\ell =0$, $\mathrm{span}_{\bc}(\{\ell_0 -\ell'\,|\,\ell'\not=\ell_0\})=\mathfrak{c}^{*}\otimes\bc$. Since $\mathrm{span}_{\bc}(N)=\mathfrak{c}^{*}\otimes\bc$, $\mathfrak{c}$ is balanced with respect to $\phi$.

From now on, we assume that $O$ is empty. If $D$ or $E$ is empty,
there is no room for pairs $(\ell,\ell')$ for $\ell\in D$ and
$\ell'\in E$, so $P$ is empty, hence balanced. Therefore we assume
that both $D$ and $E$ are non empty. Let $L_{D}$ (resp. $L_{E}$)
denote the span of all differences $\ell-\ell'$ for $\ell\in D$ and
$\ell'\in D$ (resp. for $\ell\in E$ and $\ell'\in E$). Then
$\mathrm{dim}(L_{D})=\mathrm{card}(D)-1$,
$\mathrm{dim}(L_{E})=\mathrm{card}(E)-1$ by Lemma \ref{centrslnr}.
Since $L_{D}\subset\mathrm{span}_{\bc}(D)$,
$L_{E}\subset\mathrm{span}_{\bc}(E)$ and
$\mathrm{span}_{\bc}(D)\cap\mathrm{span}_{\bc}(E)$ is the line
generated by $\sum_{\ell\in D}d_{\ell}\ell$, $L_{D}\cap
L_{E}=\{0\}$, thus $\mathrm{dim}(L_{D}+
L_{E})=\mathrm{card}(D)+\mathrm{card}(E)-2=\mathrm{dim}(\mathfrak{c})-1$.
In the quotient space $\mathfrak{c}^{*}\otimes\bc/L_{D}+ L_{E}$, all
elements of $D$ (resp. of $E$) are mapped to the same vector
$\ell_{D}$ (resp. $\ell_{E}$), and $\ell_{D}\not=\ell_{E}$. Again,
if one of the $\ell-\ell'$, $\ell\in D$, $\ell'\in E$ belongs to
$N$, $\mathfrak{c}^{*}\otimes\bc/N$ vanishes, so $\mathfrak{c}$ is
balanced with respect to $\phi$, contradiction.

Therefore, all $\ell-\ell'$, $\ell\in D$, $\ell'\in E$, belong to $\pm P$. If there exists two pairs $(\ell,\ell')$, $\ell\in D$, $\ell'\in E$, such that $\Im m(\ell-\ell')$ have opposite signs in $\mathfrak{c}^{*}/\Im m(L_{D}+ L_{E})$, then $\mathfrak{c}$ is balanced with respect to $\phi$, contradiction.

Otherwise, all $\ell-\ell'$, $\ell\in D$, $\ell'\in E$, belong to $\pm P$ and those which belong to $+P$ project to $\mathfrak{c}^{*}/\Im m(L_{D}+ L_{E})$ with equal signs. This implies that the direct sum representation $\bigoplus_{\ell\in E}I_{\ell}$ is maximal. In other words, $\phi(\Gamma)\subset U(p,p)\times U(q-p)$ is maximal. The symmetric spaces $\mathcal{D}_{p,p}$ and $\mathcal{D}_{p,q}$ corresponding to $SU(p,p)$ and $G=SU(p,q)$ have equal ranks and the embedding $\mathcal{D}_{p,p}\hookrightarrow \mathcal{D}_{p,q}$ is isometric and holomorphic. Therefore Example \ref{max2} implies that, viewed as a homomorphism $\Gamma\to G$, $\phi$ is maximal as well.
\end{pf}

\subsection{Rigidity in \texorpdfstring{$SU(p,q)$}{}}

The centralizer of $SU(p,p)$ in $SU(p,q)$, $q>p$ is $U(q-p)$ with center $\mathfrak{c}=\mathfrak{u}(1)$ generated by
$
Z=\begin{pmatrix}
-\frac{2pi}{p+q}I_{q-p}      & 0  \\
0  &  \frac{(q-p)i}{p+q}I_{2p}
\end{pmatrix}$.
There is only one nonzero pair of roots $\pm i$, giving rise to the root space $Hom_{\bc}(\bc^{q-p},\bc^{2p})$. The sesquilinear form $s_i$ on $Hom_{\bc}(\bc^{q-p},\bc^{2p})$ is, in an $SU(p,p)$-invariant manner, the direct sum of $q-p$ copies of the $U(p,p)$-invariant Hermitian form on $\bc^{2p}$. Therefore the corresponding Toledo invariant is equal to $q-p$ times the Toledo invariant obtained for $q=p+1$. In this case, the centralizer of $\mathfrak{c}$ is $U(1)\times SU(p,p)$, acting on the real root space via the standard complex representation of its second factor. The representation $\rho_i$ induced on a surface subgroup $\Gamma\subset U(1)\times SU(p,p)$ is maximal if and only if the projection of $\Gamma$ to $SU(p,p)$ is maximal. Thus a maximal surface subgroup of $SU(p,p)$ is not flexible in $SU(p,q)$, $q>p$, as is well known. Such subgroups exist (Theorem \ref{BIW03}) and are known to be automatically discrete, \cite{BGG}, \cite{BIW}.

\section{Real forms of \texorpdfstring{$O(n,\bc)$}{} and \texorpdfstring{$Sp(n,\bc)$}{}}
\label{so}

\subsection{Non pure imaginary roots}

\begin{lemma}
\label{nonimaginaryroots}
Let $G$ be a real form of $O^{\epsilon}(n,\bc)$. Let $H\subset G$ be a reductive subgroup. Let $\mathfrak{c}\subset\mathfrak{g}$ be the center of its centralizer. If one of the roots of $\mathfrak{c}$ on $\bc^n$ is not pure imaginary, then $\mathfrak{c}$ is balanced.
\end{lemma}

\begin{pf}
By contradiction. Assume that $\mathfrak{c}$ is not balanced and at least one of the roots of $\mathfrak{c}$ on $\bc^n$ is not pure imaginary.
Consider the torus $\mathfrak{c}\otimes\bc\subset \mathfrak{so}^{\epsilon}(n,\bc)$ and the set $L\cup -L$ of its nonzero roots on $\bc^n$ (Proposition \ref{centrcentrorth}). Write $L=L_I \cup L_N$ where $L_I \subset L$ denotes the set of roots whose restriction to $\mathfrak{c}$ takes only pure imaginary values and $L_N$ its complement. By assumption, $L_N$ is nonempty.

Assume first that $L_I$ is nonempty too. None of the roots $\lambda=\pm\ell\pm\ell'$, $\ell\in L_I$, $\ell'\in L_N$, is pure imaginary, so none of them belongs to $\pm P$. They generate $\mathrm{span}(L)=\mathfrak{c}^*$, thus $\mathfrak{c}$ is balanced, contradiction. Therefore, $L=L_N$.

If every root of adjoint action is non pure imaginary, it is
balanced. Hence they do not span $\mathfrak c^*$, which implies that
there exists $v\in\mathfrak c$ such that  $\lambda(v)=0$ for all
roots $\lambda$ of $\mathfrak{c}$ in the adjoint representation
which are not pure imaginary. For every distinct $\ell$ and
$\ell'\in L$, one of $\ell-\ell'$ and $\ell+\ell'$ is not pure
imaginary, thus $\ell'(v)=\pm\ell(v)$.

Assume that there exists $\ell\in L$ such that $\ell$ or $2\ell$ is
a root of an adjoint action, then $\ell(v)=0$ as well, which implies
that $\ell'(v)=0$ for all $\ell'\in L$, and $v=0$. In other words,
in that case, non pure imaginary roots span $\mathfrak{c}^*$, and
$\mathfrak{c}$ is balanced, contradiction. Therefore,
\begin{itemize}
  \item $0$ is not a root of $\mathfrak{c}$ on $\bc^n$ (i.e., all root spaces correspond to elements of $\pm L$);
  \item for every root $\ell\in L$, $2\ell$ is not a root of $\mathfrak{c}$ in the adjoint representation.
  \end{itemize}
This implies that $\epsilon=1$ and all root spaces $I_{\ell}$ have
dimension $1$, i.e. $\mathfrak{c}\otimes\bc$ is a maximal torus of
$\mathfrak{so}(n,\bc)$.  Its centralizer in $\mathfrak{so}(n,\bc)$
is abelian, so are $Z_{G}(\mathfrak{c})$ and $H$, up to finite
index. A homomorphism from an abelian group cannot be tight, so no
symplectic action on root spaces can be maximal, $P$ is empty, and
$\mathfrak{c}$ is balanced again, contradiction.
\end{pf}

\subsection{The sesquilinear structure}

\begin{prop}
\label{unitary}
Let $G$ be a real form of $O^{\epsilon}(n,\bc)$. Let $H\subset G$ be a reductive subgroup. Let $\mathfrak{c}\subset\mathfrak{g}$ be the center of its centralizer.
\begin{enumerate}
  \item Let $\ell$ be a pure imaginary root of $\mathfrak{c}$ on $\bc^n$. Let $I_{\ell}$ denote the corresponding root space. Then $I_{\ell}$ inherits a $Z_{G}(\mathfrak{c})$-invariant nondegenerate sesquilinear form $s_{\ell}$.
  \item Assume that all roots of $\mathfrak{c}$ on $\bc^n$ are pure imaginary. Then the centralizer $Z_{G}(\mathfrak{c})$ of $\mathfrak{c}$ in $G$ is a product of the fixator of the orthogonal of the $0$ root space $I_0$ and of unitary groups,
  \begin{eqnarray*}
Z_{G}(\mathfrak{c})\cong G_{I_{0}^{\bot}}\times\prod_{\ell\in L}U(I_{\ell},s_{\ell}).
\end{eqnarray*}
\end{enumerate}
\end{prop}

\begin{pf}
1. Let $s_{\ell}(v,v')=B(\tau(v),v')$. This is an $\eta\epsilon$-symmetric sesquilinear form. It is nondegenerate on $I_{\ell}$ because $B$ is and all root spaces but $I_{-\ell}=\tau(I_{\ell})$ are $B$-orthogonal to $I_{\ell}$. If $g\in Z_{G}(\mathfrak{c})$, i.e $g$ belongs to $G$ and commutes with $\mathfrak{c}$, then $g$ leaves all root spaces $I_{\pm\ell}$ invariant, it commutes with $\tau$, and is isometric for $B$, thus it is isometric for $s_{\ell}$.

2. Let $g\in G$ commute with $\mathfrak{c}$. Let $g_0 \in GL(n,\bc)$
be the element which coincides with $g$ on $I_0$ and fixes the sum
of all $I_{\ell}$, $\ell\in\pm L$, i.e. $I_{0}^{\bot}$. Then $g_0
\in G$. Indeed, $g_{0}$ commutes with $\tau$ and preserves $B$.
Also, restrict $g$ to each $I_{\ell}$, $\ell\in L$. This yields an
injective homomorphism $Z_{G}(\mathfrak{c})\to
G_{I_{0}^{\bot}}\times\prod_{\ell\in L}U(I_{\ell},s_{\ell})$.

Conversely, let $g_0 \in G$ be an element which fixes all $I_{\ell}$. For each $\ell\in L$, pick an element $g_{\ell}\in U(I_{\ell},s_{\ell})$, extend it to  $I_{-\ell}$ so that $g_{\ell}=\tau\circ g_{\ell}\circ\tau^{-1}$, extend this map trivially to other $I_{\ell'}$ and to $I_0$. The obtained linear map belongs to $Z_{G}(\mathfrak{c})$. Indeed, it preserves $B$, commutes with $\tau$ and preserves each $I_{\ell}$. Multiplying $g_0$ with the $g_{\ell}$'s yields the inverse isomorphism $G_{I_{0}^{\bot}}\times\prod_{\ell\in L}U(I_{\ell},s_{\ell})\to Z_{G}(\mathfrak{c})$.
\end{pf}

\begin{lemma}
\label{traces}
Let $G$ be a real form of $O^{\epsilon}(n,\bc)$. Let $H\subset G$ be a reductive subgroup. Let $\mathfrak{c}\subset\mathfrak{g}$ be the center of its centralizer. Let $\ell$, $\ell'$ be distinct pure imaginary roots of $\mathfrak{c}$ on $\bc^n$, and $\lambda=\ell-\ell'$ the corresponding pure imaginary root of $\mathfrak{c}$ in its adjoint action on $\mathfrak{g}$. Let $I_{\ell}$, $I_{\ell'}$ and $\mathfrak{g}_{\lambda}$ denote the corresponding root spaces.

1. If $\ell'\not=-\ell$, then, as a sesquilinear $Z_{G}(\mathfrak{c})$-module, $\mathfrak{g}_{\lambda}$ is isomorphic to $Hom(I_{\ell},I_{\ell'})$ equipped with the natural sesquilinear form
\begin{eqnarray*}
(f,f')\mapsto \mathrm{Trace}(f^{*}\circ f'),
\end{eqnarray*}
where $f^{*}\in Hom(I_{\ell'},I_{\ell})$ is the adjoint of $f$ with respect to the sesquilinear forms $s_{\ell}$ and $s_{\ell'}$.

2. If $\ell'=-\ell$, then, as a sesquilinear
$Z_{G}(\mathfrak{c})$-module, $\mathfrak{g}_{\lambda}$ is isomorphic
to the subspace of $-\epsilon$-symmetric forms in the space
$I_{\ell}^{*}\otimes I_{\ell}^{*}$ of $\bc$-bilinear forms on
$I_{\ell}$, equipped with its natural sesquilinear form
\begin{eqnarray*}
(b,b')\mapsto \mathrm{Trace}((\sharp_{s_{\ell}})^{-1}\circ(\overline{\sharp_{b}})^{\top}\circ (\sharp_{s_{\ell}})^{-1}\circ\sharp_{b'}).
\end{eqnarray*}\end{lemma}

\begin{pf}
1. Assume first that $\ell$ and $\ell'$ are linearly independant. According to Corollary \ref{rootorth}, $\mathfrak{g}_{\lambda}$ is the space of $B$-skew-symmetric elements of $Hom(I_{\ell},I_{\ell'})\oplus Hom(I_{-\ell'},I_{-\ell})$. Given $f\in Hom(I_{\ell},I_{\ell'})$ and $g\in Hom(I_{-\ell'},I_{-\ell})$, $X=(f,g)$ is $B$-skew-symmetric means that for all $v\in I_{\ell}$ and $w\in I_{-\ell'}$,
\begin{eqnarray*}
B(f(v),w)+B(v,g(w))=0.
\end{eqnarray*}
Since $\sharp_B$ identifies $I_{-\ell}$ with the dual of $I_{\ell}$, given $f\in Hom(I_{\ell},I_{\ell'})$, there exists a unique $g\in Hom(I_{-\ell'},I_{-\ell})$ such that $X=(f,g)$ is $B$-skew-symmetric. This shows that, as a $Z_{G}(\mathfrak{c})$-module, $\mathfrak{g}_{\ell-\ell'}$ is isomorphic to $Hom(I_{\ell},I_{\ell'})$.

Since $\sigma(X)=(\tau\circ g\circ\tau^{-1},\tau\circ f\circ\tau^{-1})$,
\begin{eqnarray*}
\sigma(X)\circ X=\eta(\tau\circ g\circ\tau\circ f,\tau\circ f\circ\tau\circ g),
\end{eqnarray*}
therefore
\begin{eqnarray*}
\mathrm{Trace}(\sigma(X)\circ X)&=&\eta\mathrm{Trace}(\tau\circ g\circ\tau\circ f)+\eta\mathrm{Trace}(\tau\circ f\circ\tau\circ g)\\
&=&2\eta\mathrm{Trace}(\tau\circ g\circ\tau\circ f).
\end{eqnarray*}
If $v\in I_{\ell}$,
\begin{eqnarray*}
s_{\ell}(v,\tau\circ g\circ\tau\circ f(v))
&=&B(\tau(v),\tau\circ g\circ\tau\circ f(v))\\
&=&\overline{B(v,g\circ\tau\circ f(v))}\\
&=&-\overline{B(f(v),\tau\circ f(v))}\\
&=&-\epsilon\overline{B(\tau\circ f(v),f(v))}\\
&=&-\epsilon\overline{s_{\ell}(f(v),f(v))}\\
&=&-\epsilon\overline{s_{\ell}(v,f^{*}\circ f(v))}.
\end{eqnarray*}
Summing over an orthogonal basis for $s_{\ell}$ yields
\begin{eqnarray*}
\mathrm{Trace}(\tau\circ g\circ\tau\circ
f)=-\epsilon\overline{\mathrm{Trace}(f^{*}\circ f)},
\end{eqnarray*}
hence the Killing sesquilinear form on $\mathfrak{g}_{\lambda}$ is proportional to $\mathrm{Trace}(f^{*}\circ f)$.

2. Let $f\in Hom(I_{\ell},I_{-\ell})$. Then $f$ belongs to $\mathfrak{so}(V,B)$ ($f$ is $B$-skew-symmetric) if and only if $f^{*}=-\sigma(f)$. Indeed,
\begin{eqnarray*}
&&\forall v,\,w\in I_{\ell},\quad B(f(v),w)+B(v,f(w))=0\\
&\Leftrightarrow&\forall v,\,w\in I_{\ell},\quad s_{-\ell}(f(v),\tau(w))+s_{\ell}(v,\tau\circ f(w))=0\\
&\Leftrightarrow&\forall v\in I_{\ell},\,\forall v'\in I_{-\ell},\quad s_{-\ell}(f(v),v')+s_{\ell}(v,\tau\circ f\circ \tau^{-1}(v'))=0\\
&\Leftrightarrow& f^{*}=-\tau\circ f\circ\tau^{-1}.
\end{eqnarray*}
Thus, for $B$-skew-symmetric $f$,
\begin{eqnarray*}
\mathrm{Trace}(\sigma(f)\circ f)=\mathrm{Trace}(\tau\circ f\circ\tau^{-1}\circ f)=-\mathrm{Trace}(f^{*}\circ f).
\end{eqnarray*}
Now $f$ is $B$-skew-symmetric if and only if the bilinear form $b(v,v')=B(f(v),v')$ on $I_{\ell}$ is (-$\epsilon$)-symmetric. Since $\sharp_{B}=f^{\top}\circ\sharp_{B}=\eta f^{\top}\circ\tau^{\top}\circ\sharp_{s_{\ell}}$ and $f^{*}=(\sharp_{s_{\ell}})^{-1}\circ(\overline{f})^{\top}\circ \sharp_{s_{\ell}}$,
\begin{eqnarray*}
f^{*}\circ f'=\eta(\sharp_{s_{\ell}})^{-1}\circ(\overline{\sharp_{b}})^{\top}\circ (\sharp_{s_{\ell}})^{-1}\circ\sharp_{b'}.
\end{eqnarray*}
\end{pf}

\begin{lemma}
\label{signbil} Let $V$ be an $n$-dimensional complex vectorspace
equipped with a nondegenerate sesquilinear form of signature $s$.
Then the signatures of the induced sesquilinear forms on $V^*
\otimes V^*$, $S^2 V^*$ and $\Lambda^2 V^*$ are equal to $s^2$,
$(s^2 +n)/2$ and $(s^2 -n)/2$ respectively.
\end{lemma}

\begin{pf}
Fix a basis of $V$. If $D$ denotes the matrix of the sesquilinear
form $S$ in this basis, and $b$, $b'$ the matrices of two bilinear
forms on $V$, the induced sesquilinear form on $V^* \otimes V^*$ is
\begin{eqnarray*}
\mathrm{Trace}(\bar{b}^{\top}D^{-1}b'D^{-1}).
\end{eqnarray*}
One can assume that $S$ is diagonal with entries $d_m =1$ ($p$ times) and $-1$ ($q$ times), $n=p+q$, $s=p-q$. Then
\begin{eqnarray*}
\mathrm{Trace}(\bar{b}^{\top}D^{-1}b'D^{-1})=\sum_{m,m'}d_{m}d_{m'}|b_{mm'}|^2
\end{eqnarray*}
is diagonal again. The signature of the whole space $V^* \otimes V^*$ is $\sum_{m,m'}d_{m}d_{m'}=s^2$. The signature of the subspace $\Lambda^2 V^*$ of skew-symmetric forms is
\begin{eqnarray*}
\sum_{m<m'}d_{m}d_{m'}=\frac{p(p-1)}{2}+\frac{q(q-1)}{2}-pq=\frac{s^2 -n}{2}.
\end{eqnarray*}
The signature of the subspace $S^2 V^*$ of symmetric forms is
\begin{eqnarray*}
\sum_{m\leq m'}d_{m}d_{m'}=n+\sum_{m<m'}d_{m}d_{m'}=\frac{s^2 +n}{2}.
\end{eqnarray*}
\end{pf}

\subsection{Unbalanced centers of centralizers}

Recall that, when $\mathfrak{c}\subset\mathfrak{g}$ is the center of the  centralizer of a homomorphism of a surface group to $G$, $P$ denotes the set  of roots $\lambda$ in the adjoint representation such that the sesquilinear action on the root space $\mathfrak{g}_{\lambda}$ is maximal with positive Toledo invariant.

\begin{lemma}
\label{signiell} Let $G$ be a real form of $O^{\epsilon}(n,\bc)$.
Let $\phi:\Gamma\to G$ be a reductive homomorphism. Let
$\mathfrak{c}\subset\mathfrak{g}$ be the center of its centralizer.
Assume all the roots are pure imaginary. Let $\ell$, $\ell'$ be
distinct nonzero pure imaginary roots of $\mathfrak{c}$ on $\bc^n$,
and $\lambda=\ell-\ell'$ the corresponding pure imaginary root of
$\mathfrak{c}$ in its adjoint action on $\mathfrak{g}$.
\begin{enumerate}
  \item If $\ell'\not=-\ell$, and $\lambda\in \pm P$, then one of $s_{\ell}$ and $s_{\ell'}$ has vanishing signature and the other is definite.
  \item If $\ell'=-\ell$, then $\lambda=2\ell$ either is not a root or does not belong to $\pm P$.
\end{enumerate}
\end{lemma}

\begin{pf}
If $\lambda\in P$, $\mathfrak{g}_{\lambda}$ has vanishing signature. So does $U(Hom(I_{\ell},I_{\ell'}))$ (resp. $U(\Lambda^{\epsilon}(I_{\ell}))$) with its natural sesquilinear form, according to Lemma \ref{traces}. In particular, these groups are of tube Hermitian type. Furthermore, the $Z_{G}(\mathfrak{c})$ action on $\mathfrak{g}_{\lambda}$ is tight. According to Proposition \ref{unitary}, $Z_{G}(\mathfrak{c})$ is a product of groups.

1. If $\ell'\not=-\ell$, among the factors, only $U(I_{\ell},s_{\ell})$ and $U(I_{\ell'},s_{\ell'})$ act non trivially on $\mathfrak{g}_{\lambda}$, thus the morphism
\begin{eqnarray*}
U(I_{\ell})\times U(I_{\ell'})\to U(Hom(I_{\ell},I_{\ell'}))
\end{eqnarray*}
must be tight. Lemma \ref{lemsupq} applies and one of the left hand groups is compact and the other has vanishing signature.

2. If $\ell'=-\ell$, among the factors, only $U(I_{\ell},s_{\ell})$ acts non trivially on $\mathfrak{g}_{\lambda}$, thus the morphism
\begin{eqnarray*}
U(I_{\ell})\to U(\Lambda^{\epsilon}(I_{\ell}))
\end{eqnarray*}
must be tight. Lemma \ref{tighttight} implies that $U(I_{\ell})$ is of tube type, so $s_{\ell}$ has vanishing signature. Lemma \ref{signbil} shows that the signature of $\Lambda^{\epsilon}(I_{\ell})$ is $\pm\mathrm{dim}(I_{\ell})/2$, which does not vanish, contradiction. We conclude that if $2\ell$ is a root, it does not belong to $\pm P$.
\end{pf}

\begin{prop}
\label{end} Let $G$ be a real form of $O^{\epsilon}(n,\bc)$. Let
$\phi:\Gamma\to G$ be a reductive homomorphism. Let
$\mathfrak{c}\subset\mathfrak{g}$ be the center of its centralizer.
Then $\mathfrak{c}$ is balanced with respect to $\phi$ unless
$\epsilon=1$, $\mathrm{dim}(\mathfrak{c})=1$, $\mathfrak{c}^*$ is
generated by a root $\ell$ with a $1$-dimensional root space
$I_{\ell}$, the sesquilinear form on $I_0$ has vanishing signature,
$G_{I_0^\perp}$ is reductive Hermitian of tube type and tightly
embedded in $U(I_{0})$, and the homomorphism $\Gamma\to U(I_{0})$ is
maximal.
\end{prop}

\begin{pf}
Assume that $\mathfrak{c}$ is not balanced. From Lemma \ref{nonimaginaryroots}, we know that roots have to be pure imaginary. We take their imaginary parts without mentioning it explicitely.

If $\epsilon=-1$ or if all $I_{\ell}$, $\ell\in L$, have dimension $>1$, then all $2\ell$ are roots and do not belong to $\pm P$, so they all belong to the complement $N$ of $P$ (Lemma \ref{signiell}). Since they generate $\mathfrak{c}^*$, $\mathfrak{c}$ is balanced, contradiction. So $\epsilon=1$ and the set $D$ of roots $\ell\in L$ such that $\mathrm{dim}(I_{\ell})=1$ is nonempty.

Assume that $D$ has at least $2$ elements. According to Lemma \ref{signiell}, combinations $\pm\ell\pm\ell'$, for $\ell$, $\ell'\in D$, belong to $N$. Since they span $\mathrm{span}(D)$, and multiples $2\ell''$, $\ell''\notin D$ span $\mathrm{span}(L\setminus D)$, $N$ spans $\mathfrak{c}^{*}$ and $\mathfrak{c}$ is balanced, contradiction. So $D$ has exactly one element, denoted by $\ell_0$.

Assume that $L\not=D$. For $\ell\notin D$, $2\ell\notin \pm P$,
hence $2\ell\in N$. Since $\ell_0\cup\{\ell\notin D\}$ are all
roots,  $\mathrm{dim}(\mathfrak c^{*}/\mathrm{span}(N))=1$, and for
all $\ell\in L\setminus D$, all combinations $\pm\ell_0 \pm\ell$
belong to $\pm P$. Since $0$ does not belong to the convex hull of
the image of $P$ in $c^{*}/\mathrm{span}(N)$, this convex hull
contains exactly one of $\ell_0$ and $-\ell_0$, say $\ell_0$. Then
$\ell_0 -\ell$ and $\ell_0 +\ell$ belong to $P$. Thus the
sesquilinear representation of $\Gamma$ in
$Hom(I_{\ell_0},I_{\ell}\oplus I_{-\ell})$ is maximal. Lemma
\ref{projective} allows to replace $Hom(I_{\ell_0},I_{\ell}\oplus
I_{-\ell})$ with $I_{\ell}\oplus I_{-\ell}$. However, as a
sesquilinear vectorspace, $I_{-\ell}$ is isomorphic to
$\overline{I_{-\ell}}$, so, with Lemma \ref{conjsum},
\begin{eqnarray*}
T(I_{\ell}\oplus I_{-\ell})=T(I_{\ell})+T(I_{-\ell})=T(I_{\ell})-T(I_{\ell})=0,
\end{eqnarray*}
contradicting maximality. So $L=D$ and $\mathrm{dim}(\mathfrak{c})=1$.

Assume that $0$ is not a root of $\mathfrak{c}$ on $\bc^n$. Then
$\mathrm{dim}(V=I_{\ell_0}\oplus I_{-\ell_0})=2$ and
$Z_{G}(\mathfrak{c})=U(I_{\ell_0})\cong U(1)$, which cannot have any
maximality property. So $0$ is a root on $\bc^n$, $\pm\ell_0$ are
roots of $\mathfrak{c}$ in the adjoint representation. Since $0$
does not belong to the convex hull of $P$ in $c^{*}$, exactly one of
$\pm\ell_0$ belongs to $P$, say $\ell_0 \in P$. Then
$Hom(I_{\ell_0},I_{0})$ is a maximal sesquilinear representation of
$\Gamma$, so does $I_{0}$, by Lemma \ref{projective}. In particular,
the signature of $(I_{0},s_{0})$ vanishes, i.e. $U(I_{0})$ has tube
type. Also, the morphism $Z_{G}(\mathfrak{c})\to PU(I_{0})$ is
tight. $Z_{G}(\mathfrak{c})=G_{I_0^\perp}\times U(I_{\ell_0})$ acts
on $I_{0}$ via the group $G_{I_{0}^{\bot}}$, so the injection
$G_{I_{0}^{\bot}}\to PU(I_{0})$ is tight. This implies that
$G_{I_{0}^{\bot}}$ is reductive Hermitian and has tube type (Lemma
\ref{tighttight}).
\end{pf}

\begin{cor}
\label{fin}
Let $G$ be a real form of $O^{\epsilon}(n,\bc)$. Let $\phi:\Gamma\to G$ be a reductive homomorphism. Let $\mathfrak{c}\subset\mathfrak{g}$ be the center of its centralizer. Then $\mathfrak{c}$ is balanced with respect to $\phi$ unless $G=SO^{*}(2n)$, $n$ odd, $\phi(\Gamma)\subset SO^{*}(2n-2)\times SO^{*}(2)$ and $\phi$ is maximal.
\end{cor}

\begin{pf}
In view of Proposition \ref{end}, there merely remains to determine
which pairs $(G,G_{I_{0}^{\bot}})$ can lead to unbalanced
centralizers, when $G$ is a real form of $O(N,\bc)$, i.e. $G=
O(p,q)$ or $G=SO^* (2n)$.

In the real case, $I_{0}$ is real, $G_{I_{0}^{\bot}}$ is a real
orthogonal group $O(p,q)$, which tightly injects into
$U(I_{0},s_{0})=U(p,q)$. We also know that the signature $p-q$ of
$s_0$ vanishes. $O(p,p)$ is reductive Hermitian only if $p=2$. But
Lemma \ref{o22} states that the inclusion $O(2,2)\hookrightarrow
U(2,2)$ is not tight. So $\mathfrak{c}$ is always balanced if
$G=O(p,q)$ is a real orthogonal group.

In the quaternionic case, $I_{0}$ is quaternionic and carries a
nondegenerate $(\bar{\hskip1em},-1)$-binary form (Lemma
\ref{conjugate}) (see subsection \ref{rf} also). Therefore
$G_{I_{0}^{\bot}}$ is a quaternionic skew-unitary group $SO^* (2p)$,
$p=\mathrm{dim}_{\bh}(I_{0})$. Let $\ell$ denote the unique nonzero root of $\mathfrak{c}$ on $\bc^{2n}$. $I_{0}^{\bot}=I_{\ell}\oplus I_{-\ell}$ is a $2$-dimensional complex vectorspace, thus a $1$-dimensional quaternionic vectorspace, therefore $n=p+1$. Since $G_{I_{0}^{\bot}}$ has tube type, $p$ is even and $n$ is odd. The
homomorphism $\Gamma\to U(I_{0})\cong U(p,p)$ is maximal. According
to Lemma \ref{sosu}, the homomorphism $\Gamma\to
G_{I_{0}^{\bot}}\cong SO^{*}(2p)$ is maximal as well. The symmetric
spaces $\mathcal{G}_{p}$ and $\mathcal{G}_{p+1}$ corresponding to
$G_{I_{0}^{\bot}}$ and $G$ have equal ranks and the embedding
$\mathcal{G}_{p}\hookrightarrow \mathcal{G}_{p+1}$ is isometric and
holomorphic. Therefore Example \ref{max3} implies that, viewed as a
homomorphism $\Gamma\to G$, $\phi$ is maximal as well.
\end{pf}

\section{Proof of Theorem \ref{classical}}
\label{nonreductive}

For homomorphisms $\Gamma\to G$ with reductive Zariski closure, the proof of Theorem \ref{classical} follows from Theorem \ref{flexssimple}, the classification of classical simple Lie groups and the case by case analysis of balancedness in sections \ref{sl} and \ref{so}.

Here is how the problem is reduced to the case of reductive homomorphisms. Theorem 2 of \cite{KP} asserts that if $\mathrm{genus}(\Gamma)\geq 2\mathrm{dim}(G)^2$ and $G$ is semisimple, the space $Hom(\Gamma,G)$ falls into two types of connected components: in some of them, Zariski dense homomorphisms are dense; others do not contain any Zariski dense homomorphism (call them rigid). Let $\phi:\Gamma\to G$ be an arbitrary homomorphism. If $\phi$ cannot be approximated by Zariski dense homomorphisms, then $\phi$ belongs to a rigid component $\mathcal{C}$. Proposition 8.3 of \cite{KP} asserts that $\mathcal{C}$ contains a reductive homomorphism $\psi$. The reductive case of Theorem \ref{classical} implies that $G$ is Hermitian of non tube type and $\psi$ is maximal. Since Toledo invariants are constant on connected components of $Hom(\Gamma,G)$, $\phi$ is maximal as well.

\section{Appendix}
\label{app}

For the reader's convenience, we give barehanded proofs of the maximality preserving property of two embeddings between reductive Hermitian groups.

We start with a preliminary observation.

\begin{lemma}
\label{poly}
Let $F:Y\to X$ be an equivariant totally geodesic map between Hermitian symmetric spaces. Assume $Y$ is irreducible. Let $P\subset X$ and $Q\subset Y$ be maximal polydisks such that $F(Q)\subset P$. Then $F$ is positively maximality preserving if and only $F_{|Q}:Q\to P$ is.
\end{lemma}

\begin{pf}
Since $Y$ is irreducible, $F$ is homothetic, i.e. there exists a constant $c$ such that $F^{*}\omega_X =c\,\omega_{Y}$. Sectional curvature achieves its minimum along maximal polydiscs, so
\begin{eqnarray*}
(\omega_{X})_{|P}=\omega_{P},\quad (\omega_{Y})_{|Q}=\omega_{Q}.
\end{eqnarray*}
Since $\mathrm{rank}(P)=\mathrm{rank}(X)$ and $\mathrm{rank}(Q)=\mathrm{rank}(Y)$, $F$ is positively maximality preserving if and only if
\begin{eqnarray*}
c=\frac{\mathrm{rank}(X)}{\mathrm{rank}(Y)}
&\Leftrightarrow&
c=\frac{\mathrm{rank}(P)}{\mathrm{rank}(Q)}\\
&\Leftrightarrow&
(F_{|Q})^{*}\omega_{P}=\frac{\mathrm{rank}(P)}{\mathrm{rank}(Q)}\omega_{Q}
\end{eqnarray*}
if and only if $F_{|Q}$ is positively maximality preserving.
\end{pf}

\begin{lemma}
\label{max1}
The embedding $SU(n,n)\hookrightarrow Sp(4n,\br)$, is positively maximality preserving.
\end{lemma}

\begin{pf}
Let $Y=\mathcal{D}_{n,n}$ (resp. $X=\mathcal{S}_{2n}$) denote the
symmetric space associated to $H=SU(n,n)$  (resp. $G=Sp(4n,\br)$).
Let $\iota:Y\to X$ denote the corresponding embedding of symmetric
spaces. We must show that
\begin{eqnarray*}
\iota^{*}\omega_{X}=2\omega_{Y}.
\end{eqnarray*}

Let us first study the case when $n=1$. Let $V_{\bc}=\bc^{2}$ be
equipped with the standard symmetric  sesquilinear form $v\cdot
v'=\bar{v}^{\top}v'$. Let $S=\begin{pmatrix}
0 & i\\
-i & 0
\end{pmatrix}$. The symmetric sesquilinear form $s(v,v')=v\cdot (Sv')$ on $V_{\bc}$ has vanishing
signature. It is easy to show that $s(Av, Av')=s(v,v')$ for
$v,v'\in\bc^2$, and $A\in SL(2,\br)$ by a direct calculation. Hence
the group $H=SU(V_{\bc},s)$ coincides with $SL(2,\br)$ acting on
$V_{\bc}=\br^2 \otimes\bc$. Its maximal compact subgroup $L$ is
generated by $J=\begin{pmatrix}
0 & 1\\
-1 & 0
\end{pmatrix}$.

Let $V_{\br}$ denote $V_{\bc}$ viewed as a real vectorspace equipped with the symplectic form $\Omega(v,v')=\Im m(s(v,v'))$. Then $H=SU(V_{\bc},s)$ is a subgroup of the larger symplectic group $G=Sp(V_{\br},\Omega)$. Let $\rho:SU(V_{\bc},s)\hookrightarrow Sp(V_{\br},\Omega)$ denote the inclusion homomorphism. $J'=\rho(J)$ is the $4\times 4$ matrix which reads $\begin{pmatrix}
0 & 1\\
-1 & 0
\end{pmatrix}$ in $2\times 2$ blocks. $J'$ is a complex structure compatible with $\Omega$ and tamed by $\Omega$, thus its centralizer in $Sp(V_{\br},\Omega)$ is a maximal compact subgroup $K$ of $Sp(V_{\br},\Omega)$. The adjoint actions of $J$ on $\mathfrak{h}/\mathfrak{l}$ and of $J'$ on $\mathfrak{g}/\mathfrak{k}$ define the complex structures of the symmetric spaces $Y$ and $X$ associated to $H$ and $G$, so the $\rho$-equivariant embedding $\iota:Y\hookrightarrow X$ mapping $L$ into $K$ is holomorphic.

Let us view $V_{\bc}$ as $\br^2 \otimes\bc$. Then $V_{\br}=\br^2
\oplus i\br^2$. In this coordinates, $SL(2,\br)$ acts diagonally.
The stabilizer of this decomposition in $G$ is the standard
$SL(2,\br)\times SL(2,\br)$ of $Sp(4,\br)$. By Example \ref{max4},
the embedding of symmetric spaces $Y=\mathcal{S}_{1}\hookrightarrow
\mathcal{S}_{2}$ corresponding to each $SL(2,\br)\hookrightarrow G$
is isometric and holomorphic, so $SL(2,\br)\times
SL(2,\br)\hookrightarrow G$ gives rise to an isometric and
holomorphic map of $Y\times Y$ onto a maximal polydisk $P$ of
$X=\mathcal{S}_{2}$. The image $\rho(H)$ sits diagonally in the
product $SL(2,\br)\times SL(2,\br)$, so $\iota:Y\to P$ factors
through the diagonal $\Delta:Y\to Y\times Y$. This shows that
K\"ahler forms fit up to a factor of $2$, i.e.
\begin{eqnarray*}
\iota^{*}\omega_{P}=\Delta^{*}(pr_{1}^{*}\omega_{P}+pr_{2}^{*}\omega_{P})=2\omega_{Y}.
\end{eqnarray*}

In general, let $(V_{\bc},s)$ be the orthogonal direct sum of $n$
copies of the $n=1$ example just studied.  Then $s$ has vanishing
signature. Each factor gives rise to a homomorphism
$SU(1,1)\hookrightarrow H=SU(V_{\bc},s)$ and a map
$\mathcal{D}_{1,1}=\mathcal{S}_{1}\hookrightarrow
Y=\mathcal{D}_{n,n}$ which, according to Example \ref{max2}, is
isometric and holomorphic. The product map
$\mathcal{D}_{1,1}^{n}\hookrightarrow Y$ is isometric and
holomorphic onto a maximal polydisk $Q$ of $Y$.

Let $V_{\br}$ be $V_{\bc}$ viewed as a real vectorspace equipped
with the symplectic structure $\Omega=\Im m(s)$.  Each factor of
$V_{\bc}$ is the complexification of a real $2$-dimensional
vectorspace. This gives rise to commuting embeddings
$SL(2,\br)\hookrightarrow Sp(4n,\br)$ and the corresponding map
$\mathcal{S}_{1}^{2n}\hookrightarrow X=\mathcal{S}_{2n}$ is an
holomorphic isometry onto a maximal polydisk $P$ of $X$. The
restriction of $\iota:Y\to X$ to $Q$ is the direct product of $n$
copies of the $n=1$ case, so again
\begin{eqnarray*}
\iota^{*}\omega_{P}=2\omega_{Q}.
\end{eqnarray*}
With Lemma \ref{poly}, since $\mathrm{rank}(X)=2\mathrm{rank}(Y)$, this shows that $\iota$ is positively maximality preserving.
\end{pf}

\begin{lemma}
\label{sosu}
Let $h$ be a nondegenerate ($\bar{\hskip1em}$,skew)-symmetric binary form on a $2n$-dimensional quaternionic vectorspace $V_{\bh}$ (see subsection \ref{defcla}). Use right multiplication by $i$ to turn $V_{\bh}$ into a complex $4n$-dimensional vectorspace $V_{\bc}$. Let $\mathcal{C}(q)=a$ denote the complex part of a quaternion $q=a+jb$. Then
\begin{eqnarray*}
s(v,v')=\mathcal{C}(h(v,v'))
\end{eqnarray*}
is a nondegenerate sesquilinear form of vanishing signature on
$V_{\bc}$. Then the corresponding embedding of groups
$\rho:SO^{*}(4n):=SU(V_{\bh},h)\hookrightarrow SU(2n,2n)$
$:=SU(V_{\bc})$ is positively maximality preserving.
\end{lemma}

\begin{pf}
Let $\iota:Y=\mathcal{G}_{2n}\to X=\mathcal{D}_{2n,2n}$ be the
corresponding embedding of Hermitian symmetric spaces. Since
$\mathrm{rank}(\mathcal{G}_{2n})=n$ and
$\mathrm{rank}(\mathcal{D}_{2n,2n})=2n$, we shall show that K\"ahler
forms match up to a factor of $2$, i.e.
$\iota^{*}\omega_{X}=2\omega_{Y}$. Following Lemma \ref{poly}, it
suffices to understand the restriction of $\iota$ to maximal
polydisks.

Let $v\cdot v'=\bar{v}^{\top}v'$ denote the standard positive definite $\bar{\hskip1em}$-symmetric binary form on $\bh^{2n}$. Let $h(v,v')=v\cdot iv'$. Then $h$ is ($\bar{\hskip1em}$,skew)-symmetric and nondegenerate, so we can take $V_{\bh}=(\bh^{2n},h)$. The embedding $\rho:SO^{*}(4n)\hookrightarrow SU(2n,2n)$ consists in taking a quaternionic matrix $X$, splitting it as $X=M+jM'$ where $M$ and $M'$ have complex entries, letting $X$ act on the quaternionic vector $v=a+jb$, $(a,b)\in(\bc^{2n})^2 =\bc^{4n}$. Thus
\begin{eqnarray*}
Xv=(M+jM')(a+jb)=Ma-\overline{M'}b+j(M'a+\bar{M}b),
\end{eqnarray*}
i.e.
\begin{eqnarray}
\label{rho}
\rho(X)=\begin{pmatrix}
M  & -\overline{M'} \\
M' & \bar{M}
\end{pmatrix}.
\end{eqnarray}
Let $J\in Gl(2n,\bh)$ denote left multiplication by $i$. Elements of $Sp(2n)$ which commute with $J$ (i.e. matrices with entries in $\bc\subset\bh$) form a group $L$ isomorphic to $U(2n)$. It is a maximal compact subgroup in $SO^{*}(4n)$. Under $\rho$, this subgroup is mapped to the maximal compact subgroup $K=S(U(2n)\times U(2n))$ by $M\mapsto (M,\bar{M})$. $J$ belongs to the Lie algebra $\mathfrak{h}=\mathfrak{so}^{*}(4n)$, it generates the center of the Lie algebra $\mathfrak{l}$ of $L$. Therefore the complex structure on $\mathcal{G}_{2n}$ arises from the adjoint action of $J$ on $\mathfrak{h}/\mathfrak{l}$. For the same reason, the complex structure on $\mathcal{D}_{2n,2n}$ arises from the adjoint action of $J'=diag(i,\ldots,i,-i,\ldots,-i)$ on $\mathfrak{g}/\mathfrak{k}$. Note that at the Lie algebra level $\rho(J)=J'$, thus $\iota$ is holomorphic.

Let us first study the case when $n=1$. The Lie algebra
$$\mathfrak{so}^{*}(4)=\{A| A^*\begin{pmatrix}
                                i & 0 \\
                                0 & i\end{pmatrix}+\begin{pmatrix}
                                                    i & 0\\
                                                    0 & i \end{pmatrix}A=0\}$$ is isomorphic to $\mathfrak{su}(1,1)\oplus
\mathfrak{su}(2)$, it consists of matrices of the form
\begin{eqnarray*}
\begin{pmatrix}
i\alpha & jx\\
-jx & i\alpha
\end{pmatrix}
+\begin{pmatrix}
i\beta & y\\
-\bar{y} & -i\beta
\end{pmatrix}
\end{eqnarray*}
where $\alpha$, $\beta\in\br$, $x$, $y\in\bc$. The first matrix belongs to a subalgebra $\mathfrak{q}$ isomorphic to $\mathfrak{su}(1,1)$. A computation based on formula (\ref{rho}) gives
\begin{eqnarray*}
\rho\begin{pmatrix}
i\alpha & jx\\
-jx & i\alpha
\end{pmatrix}=\begin{pmatrix}
i\alpha&0&0 & -\bar{x}\\
0& i\alpha&\bar{x}&0\\
0&x&-i\alpha&0\\
-x&0&0&-i\alpha
\end{pmatrix}.
\end{eqnarray*}
We see that $\rho(\mathfrak{q})$ is contained in the subalgebra
\begin{eqnarray*}
\mathfrak{p}=\{\begin{pmatrix}
i\alpha&0&0 & -\bar{x}\\
0& i\alpha'&-\overline{x'}&0\\
0&-x'&-i\alpha'&0\\
-x&0&0&-i\alpha
\end{pmatrix}\,;\,
\alpha,\,\alpha'\in\br,\,x,\,x'\in\bc\}
\end{eqnarray*}
which is isomorphic to $\mathfrak{su}(1,1)\oplus\mathfrak{su}(1,1)$ and embedded in the standard (block diagonal) manner in $\mathfrak{su}(2,2)$. The map $\rho_{|\mathfrak{q}}:\mathfrak{q}\to\mathfrak{p}$ is the graph of an inner automorphism of $\mathfrak{su}(1,1)$. Geometrically, this means that $\iota$ maps $Q=\mathcal{G}_{2}$ (a complex line of constant curvature $-1$) holomorphically into a maximal polydisk $P$ of $\mathcal{D}_{2,2}$. $P$ is holomorphically isometric to $Q\times Q$ and $\iota_{|Q}:Q\to P$ is the graph of an isometry $I$ of $Q$. It follows that
\begin{eqnarray*}
\iota^{*}\omega_{P}=(id,I)^{*}(pr_{1}^{*}\omega_{Q}+pr_{2}^{*}\omega_{Q})=2\omega_{Q}.
\end{eqnarray*}

Let us map $SO^{*}(4)$ as a diagonal $2\times 2$ block in $SO^{*}(4n)$. This yields an embedding $\mathcal{G}_{2}\hookrightarrow \mathcal{G}_{2n}=Y$ which, according to Example \ref{max3}, is isometric and holomorphic. Splitting $V_{\bh}=\bh^{2n}$ as an orthogonal direct sum of $2$-dimensional quaternionic vectorspaces yields an isometric and holomorphic embedding of $\mathcal{G}_{2}^{n}$ onto a maximal polydisk $Q\subset\mathcal{G}_{2n}$. Then $\iota(Q)$ is contained in the standard maximal polydisk $P\subset\mathcal{D}_{2n,2n}=X$, and $\iota_{|Q}:Q\to P$ is a product of $n$ copies of the map of the previous paragraph. Therefore
\begin{eqnarray*}
\iota^{*}\omega_{P}=\sum_{i=1}^{n}\iota^{*}\omega_{P_{i}}=\sum_{i=1}^{n}\omega_{Q_{i}}
=2\omega_{Q}.
\end{eqnarray*}
Since $\mathrm{rank}(X)=2\mathrm{rank}(Y)$, this shows that $\iota$ (and thus $\rho$) is positively maximality preserving.
\end{pf}

2000 {\sl{Mathematics Subject Classification.}} 51M10, 57S25.

{\sl{Key words and phrases.}} Algebraic group, symmetric space, rigidity, group cohomology, moduli space, tube type

\vskip1cm

\noindent     Inkang Kim\\
     School of Mathematics\\
     KIAS, Heogiro 87, Dongdaemen-gu\\
     Seoul, 130-722, Korea\\
     \texttt{inkang\char`\@ kias.re.kr}

\smallskip

     \noindent  Pierre Pansu\\ D\'epartement de Math{\'e}matiques
et applications \\
     UMR 8553 du CNRS\\
 \'Ecole Normale Sup\'erieure\\
 45 rue d'Ulm\\
 75230 Paris C\'edex 05, France\\
  \texttt{pierre.pansu\char`\@ ens.fr}

     \end{document}